\documentclass[12pt]{article}
\usepackage{amsmath,epsfig,amsfonts,amssymb}
\usepackage{amssymb,amsthm,color,graphicx}
\usepackage{cancel}
\usepackage{epstopdf}

\def\ebe{\begin{equation}}
\def\een{\end{equation}}
\def\bea{\begin{eqnarray}}
\def\eea{\end{eqnarray}}

\def\etal{\mbox{\it et al.}}

\newcommand{\half}{\frac{1}{2}}

\newcommand{\nb}{\mathbf{n}}

\newcommand{\gradb}{\boldsymbol{\nabla}}

\begin{document}
\bibliographystyle{plain} 
\title{A Conservative Finite Difference Scheme for Poisson-Nernst-Planck Equations}
\author{Allen Flavell\\ \textit{\small Department of Applied Mathematics, Illinois Institute of Technology}\and Michael Machen\\ \textit{\small Department of Applied Mathematics, Illinois Institute of Technology}\and Bob Eisenberg\\ \textit{\small Department of Molecular Biophysics and Physiology, Rush University}\and Chun Liu\\ \textit{\small Department of Mathematics, Pennsylvania State University} \and Xiaofan Li\\ \textit{\small Department of Applied Mathematics, Illinois Institute of Technology}}
\maketitle

\begin{abstract}
A macroscopic model to describe the dynamics of ion transport in ion channels 
is the Poisson-Nernst-Planck(PNP) equations. In this paper, we develop 
a finite-difference method for solving PNP equations, 
which is second-order accurate in both space and time. 
We use the physical parameters specifically suited 
toward the modelling of ion channels.
We present a simple iterative scheme to solve the system of nonlinear 
equations resulting from discretizing the equations implicitly in time,
which is demonstrated to converge in a few iterations. 
We place emphasis on ensuring numerical methods to have 
the same physical properties that the PNP equations themselves also possess, 
namely conservation of total ions and correct rates of energy dissipation. 
We describe in detail an approach to derive a finite-difference method
that preserves the total concentration of ions exactly in time. Further, 
we illustrate that, using realistic values of the physical parameters,
the conservation property is critical in obtaining correct numerical solutions
over long time scales.  
\end{abstract}


\section{Introduction}
The Poisson-Nernst-Planck(PNP) equations describe the diffusion of ions under the effect of an electric field that is itself caused by those same ions. The system is created by coupling the Nernst-Planck equation (which describes the diffusion of ions under the effect of an electric potential) with the Poisson equation (which relates charge density with electric potential). This system of equations has found much use in the modelling of semiconductors.\cite{markowich1990semiconductor}  
Although the Poisson-Nernst-Planck equations were applied to model membrane transport for longer than they have been employed to model semiconductors\cite{Teorell}, the use of the system to model the behavior of the internal mechanics of these transport processes is much more recent.\cite{eisenberg1998}

The system of PNP equations and its related models have been the subject of much study and numerical simulation. A recent advancement in this field was the application of energy variational analysis and density functional theory to modify the PNP system to accommodate various phenomena exhibited by 
biological ion channels. See \cite{Wei2012} and the references therein.

The computer simulations of the Poisson-Nernst-Planck models are able 
to capture the transient, dynamical behavior of the system, 
and the numerical schemes employed are quite varied. 
Cagni \etal\ (2007) \cite{Cagni2007} discretized 
the PNP in two dimensions using a second-order accurate 
finite difference method with central differencing in space and 
Crank-Nicolson scheme in time, and simulated an ion channel 
subjected to time-dependent perturbations. 
Nanninga (2008) \cite{Nanninga2008} studied a nerve impulse using a similar
finite difference scheme as in \cite{Cagni2007} but in three dimensions, 
notable in that it directly included gating and selectivity into the model. 
Lopreore \etal\ (2008) \cite{Lopreore2008} developed a finite-volume-based 
technique to solve PNP in three dimensions, which decomposes the domain 
using a dual Delaunay-Voronoi mesh. 
Neuen (2010) \cite{Neuenthesis2010} developed a semi-implicit 
finite element-based scheme to simulate three-dimensional, multi-scale 
extended PNP. Gardner and Jones (2011) \cite{Gardner2011} 
simulated a potassium channel modelled with PNP in two dimensions 
using a finite difference method with TR-BDF2 time integration. 
Much of the numerical schemes in \cite{Gardner2011} is based 
on the previous work \cite{Gardner04electrodiffusionmodel}, 
a one-dimensional model of the same channel. 
Hyon \etal\ (2011) \cite{liu10a} presented another finite element method
with back-Euler method in time to investigate the effects of finite size
of the ions by modifying the PNP via introducing a repulsive potential energy 
into the total energy. Horng \etal\ (2012) \cite{Horng2012} applied 
the multiblock Chebyshev pseudospectral method and the method of lines to solve 
a one-dimensional modified PNP modelling the finite-sizeness of the ions via
a local model.

One of the characteristics of the nonlinear PNP equations 
is that its overall behavior is very sensitive to the boundary conditions.\cite{gillespie2002}  
This presents a challenge for accurate and efficient numerical simulations, 
as generally the boundary conditions will have to be 
discretized and approximated.
In this paper, we shall investigate the effects of discretization error 
on the Poisson-Nernst-Planck equations, in particular 
discretization of the boundary conditions and the equations at the boundaries. 
We will demonstrate that the conservation properties of the numerical
methods could be critical in obtaining the long-time behavior of the solutions.

The paper is organized as follows.
We start by defining and simplifying the equations we are working with, in Sec.~\ref{secGoverningEquations}, including the introduction of the quantities 
that shall be preserved by our numerical schemes:  
the total concentration of each ion species in Sec.~\ref{secTotalConcentration}
and the energy dissipation law in Sec.~\ref{secEnergyDissipation}. 
We then describe our numerical schemes in Sec.~\ref{secNumericalMethods}, 
which presents an approach to conserve the total ion concentrations exactly 
and approximate the energy dissipation law closely. 
Finally, we shall discuss the results of simulating the system using our numerical schemes in Sec.~\ref{secNumericalResults}.


\section{Governing Equations}
\label{secGoverningEquations}
   Consider the PNP equations \cite{eisenberg1998,Gardner04electrodiffusionmodel}
\bea
   \frac{\partial c_{i}}{\partial t} &=& \gradb \cdot 
       \left\{ D_{i} \left[ \gradb c_{i}+ \frac{z_i e}{k_B T} 
          c_i\gradb \phi \right] \right\}, i=1,2,\ldots, N,
      \label{eq.id} \\
    \gradb \cdot (\epsilon \gradb \phi) &=& - \left( \rho_0 + \sum_{i=1}^N z_i e c_{i}\right),
      \label{eq.ep} 
\eea
where $c_i$ is the ion density for the $i$-th species, $D_i$ is the diffusion constant, $z_i$ is the valence, $e$ is the unit charge, $k_B$ is the Boltzmann constant, $T$ is the absolute temperature, $\epsilon$ is the permittivity, $\phi$ is the electrostatic potential, $\rho_0$ is the permanent (fixed) charge density of the system, and $N$ is the number of ion species.\cite{liu10a} 
The equations are valid in a bounded domain $\Omega$ with boundary $\partial\Omega$
and for time $t\geq 0$. 

In this work, we shall use the no-flux boundary condition for Eq.~(\ref{eq.id}). This may correspond to modelling the interior conditions of a channel that is in an occluded state, with closed gates at either end. Simulations of channels such as the KirBac1.1 channel in such a state have been conducted in the past\cite{domene2008occluded}. We shall use the Robin boundary condition for the Poisson equation,
 which models the effects of making the source of the potential across the channel
partially removed from the ends of the channel. 
The formula for the boundary conditions are 
\begin{subequations}
 \ebe
   D_{i} \left[ \gradb c_{i} + \frac{z_i e}{k_B T}
          \gradb c_i \phi 
       \right] \cdot \nb = 0, \quad i=1,2,\ldots, N, 
\label{eq.npbc}
 \een
\ebe
 \quad (\phi - \phi_\pm) + \eta \frac{\partial \phi}{\partial \nb} = 0, 
\label{eq.pbc}
\een
\end{subequations}
for points on the boundary $x\in \partial\Omega$.

For some situations, such as a generic potassium channel separating potassium and chloride ion baths, the experimental data can be well-approximated by a one-dimensional model.\cite{Gardner04electrodiffusionmodel} In one dimension, the equations \eqref{eq.id} and \eqref{eq.ep} are simplified as

\bea
 \frac{\partial c_i}{\partial t} &=& \frac{\partial}{\partial x} 
       \left[  D_i \left( \frac{\partial c_i}{\partial x} 
          + \frac{z_i e}{k_B T} c_i \frac{\partial\phi}{\partial x}\right)
       \right]  \label{eq.1did} \\
    \frac{\partial}{\partial x} \left( \epsilon 
    \frac{\partial\phi}{\partial x} \right) 
      &=& - \left(\rho_0 +\sum_i z_i e c_i\right),
      \label{eq.1dep} 
\eea
for $-L\leq x\leq L$ and $t\geq 0$, where $L$ is the half of the length
of the ion channel.
The corresponding boundary conditions are
\ebe
    \frac{\partial}{\partial x} 
       \left[ D_i \left( \frac{\partial c_i}{\partial x} 
          + \frac{z_i e}{k_B T} c_i \frac{\partial\phi}{\partial x}\right)
       \right] = 0, \quad
   (\phi - \phi_\pm) \pm \eta \frac{\partial \phi}{\partial x} = 0, 
   \quad \text{for } x=-L,L.
\label{eq.1dbc}
\een

\subsection{Total Concentration}
\label{secTotalConcentration}
The total concentration per ion species is given by
\ebe
c_{i,tot}(t) = \int_{-L}^L c_{i}(x,t) \, {\rm d} x, \quad i=1,2,\ldots , N.
\label{eq.ctot}
\een

Due to the no-flux boundary conditions \eqref{eq.1dbc}, 
the total concentration of each ion species is constant in time. 
This can be verified easily by differentiating \eqref{eq.ctot} 
with respect to time, then applying the convection-diffusion equation \eqref{eq.1did} and no flux boundary condition \eqref{eq.1dbc}.

One of the metrics we can use to evaluate different numerical schemes is therefore to measure how well the total concentration is conserved in numerical simulation. Ensuring that the total concentration for each species $c_{i,tot}$ 
is constant will be the idea behind the schemes presented in this work. 
As will be seen in Sec.~\ref{secNumericalResults}, the preservation of the conservation property is crucial for producing correct numerical results over long time scales.


\subsection{Energy Dissipation}
\label{secEnergyDissipation}
The governing equations \eqref{eq.1did} and \eqref{eq.1dep} 
for the transport of ions 
can be derived from the energy of the system using
variational principles. Similar to \cite{liu10a}, the 
total energy for our specific system is defined by 
\begin{equation}
  E = \int_{-L}^L \left[ k_B T \sum_{i=1}^{N} c_i \log \frac{c_i}{c_{i,0}} 
  + \frac{1}{2} (\rho_0 + \sum_{i=1}^{N} z_i e c_i) \phi \right] \, {\rm d} x
  + \frac{\epsilon}{2\eta}
       ( \phi_+ \phi(L) + \phi_- \phi(-L) ), 
\label{eq.ten}
\end{equation}
where $c_{i,0}$ are constants called ``reference concentrations''. 
Using the Poisson equation \eqref{eq.1dep}, the total energy can be written as
\begin{equation}
    E = \int_{-L}^L \left[ k_B T \sum_{i=1}^{N} c_i \log \frac{c_i}{c_{i,0}} 
  + \frac{\epsilon}{2} \left( \frac{\partial \phi}{\partial x} \right)^2 \right] \, {\rm d} x
  + \frac{\epsilon}{2\eta} ( \phi^2(L) +  \phi^2(-L) ),
\label{eq.ten1}
\end{equation}
where the last term is the contribution of the electric energy from the boundaries.
The total energy $E$ satisfies the energy dissipation property
\begin{equation}
  \frac{dE}{dt}=-\int_{-L}^L \sum_{i =1}^{N} \frac{D_i}{k_B T} c_i \left| 
     \frac{\partial \mu_i}{\partial x}\right|^2 \, {\rm d} x,  
\label{eq.ourenergylaw}
\end{equation}
where $\mu_i$ is the chemical potential of $i$'th ion species defined by
the variational derivative of the energy with respect to the concentration
$c_i$
\begin{equation}
  \mu_i = \frac{\delta E}{\delta c_i} = k_B T \left(\log \frac{c_i}{c_{i,0}} + 1\right) + z_i e \phi. 
\label{eq.pi}
\end{equation}
The energy dissipation law \eqref{eq.ourenergylaw}
can be derived by taking the time derivative of the total energy~\eqref{eq.ten} and 
applying integration by parts, Eqs.~\eqref{eq.1did}-\eqref{eq.1dep} 
and the boundary condition \eqref{eq.1dbc}:
\bea 
  \frac{{\rm d}E}{{\rm d}t} &=&  \int_{-L}^L \left[ k_B T \sum_i 
    ( \log \frac{c_i}{c_{i,0}} + 1) \frac{\partial c_i}{\partial t} 
    + \frac{1}{2} \sum_i z_i e \frac{\partial c_i}{\partial t} \phi
    + \frac{1}{2} (\rho_0 + \sum_i z_i e c_i) \frac{\partial \phi}{\partial t}
    \right] \, {\rm d} x  \nonumber \\ 
 &&+ \frac{\partial}{\partial t} 
       \left[ \frac{\epsilon}{2\eta} 
       ( \phi_+ \phi(L) + \phi_- \phi(-L) ) \right] \nonumber \\
 &=& - \int_{-L}^L \sum_i \frac{D_i}{k_B T} c_i \left| 
     \frac{\partial \mu_i}{\partial x}\right|^2 \, {\rm d} x  
    - \frac{1}{2} \epsilon \left(\frac{\partial \phi}{\partial x} 
    \frac{\partial \phi}{\partial t}  
    - \frac{\partial^2 \phi}{\partial x\partial t}\phi
    \right) \bigg|_{-L}^{L}\nonumber \\
 &&+\frac{\partial}{\partial t} 
       \left[ \frac{\epsilon}{2\eta} 
       ( \phi_+ \phi(L) + \phi_- \phi(-L) ) \right]. \label{eq.ed1}
\eea
The rate of energy decay \eqref{eq.ourenergylaw} can be obtained 
by using the boundary condition \eqref{eq.1dbc} 
to show the last two terms on the RHS of \eqref{eq.ed1} cancel each other. 


\subsection{Parameters and Nondimensionalization}
\label{secNonDimensionalization}

We specify the units and the parameters using 
the approximate values corresponding to the KcsA potassium channel\cite{doyle1998}. 
In our 1D model, the cylindrical channel takes a diameter of $10$ \r{A} and a length of $120$ \r{A}. We shall assume no permanent charges or selectivity for the purposes of this simulation. We consider the case of two ion species, i.e. $N=2$, with 
the initial concentration for each ion being 2 molar, resulting in an initial
number density (number of ions per unit volume) of
$1.2044 \times 10^{-3}\ \text{ions/\r{A}}^3$. 
The combination of the parameters $k_B T/e$ is approximately
$0.025\ \text{V}$, assuming the temperature is $T=298\ \text{K}$.
The permittivity $\epsilon = \epsilon_r \epsilon_0$ 
is determined by the value of the vacuum 
$\epsilon_0=8.854187817 \times 10^{-12}\ \text{F/m}$ and the relative
permittivity $\epsilon_r$ ($78.5$ for water).

The values of the diffusion coefficients $D_i$ depend on both the ion species and the channel. The only net effect of different diffusion constants is the rate of evolution of the system. Typical values for the diffusion coefficients for ion species in a channel are around $10^{9}$ \r{A}$^2$/s.\cite{dirkdiffusion2008} We will select both diffusion coefficients to be equal to each other, causing them to take a value of one after nondimensionalization.

The parameter $\eta$, as a component of the Robin boundary condition \eqref{eq.pbc},
is an aggregate of multiple physical constants and is highly dependent on the properties of the surrounding membrane. Modelling the experimental setup as an electrical circuit shows that the quantity $A\epsilon_l/\eta$, where $A$ is the area of the membrane and $\epsilon_l$ is the permittivity of the membrane, has units of capacitance and is related to charge storage. The most significant charge storage contributing to $A\epsilon_l/\eta$ is in fact the membrane capacitance, so we may surmise that the primary contributor to $\eta$ is the membrane capacitance. If a very high capacitance to ground is present, $\eta$ is approximated by the appealing formula $\eta=A\epsilon_l/C$, where $C$ is the capacitance of the membrane, however realistically $\eta$ is much smaller than that. In this work, we shall take $\eta = 2.78\times 10^{-3}\ \text{\r{A}}$ for our numerical simulations, but will also examine the effects of $\eta$ over a range from $10^{-5} \ \text{\r{A}}$ to $60\ \text{\r{A}}$. Changing the value of $\eta$ might correspond to adding a parallel capacitance in experiment. 


  Define the dimensionless variables and parameters $c'_i=c_i/c_0$, $x'=x/L$, 
$t'=t/(L^2/D_0)$, $D'_i=D_i / D_0$, $\phi'=\phi/\phi_0$,
where $c_0$ is the average of the initial charge concentration,
$L$ is the half of the channel length or computational domain, $D_0$ is 
a typical diffusion coefficient, $\phi_0$ is a characteristic value
of the electrostatic potential such as the boundary value. 
Then, non-dimensionalizing the Nernst-Planck Eq.~\eqref{eq.1did}, we obtain 
\ebe
   \frac{\partial c'_i}{\partial t'} = \frac{\partial}{\partial x'} 
       \left\{ D'_i \left[ \frac{\partial c'_i}{\partial x'} 
          + \chi_1 \left( z_i c'_i \frac{\partial\phi'}{\partial x'} \right)
       \right] \right\}, \quad \text{where } \chi_1 := e\phi_0/k_BT.
\label{eq.nid} 
\een
From the above, the dimensionless parameter $\chi_1 \approx 3.1$, 
if $\phi_0=0.08$V. 
The nondimensionalized Poisson Eq.~\eqref{eq.1dep} is given by
\ebe
    \frac{\partial}{\partial x'} \left(\epsilon' 
    \frac{\partial\phi'}{\partial x'} \right) 
      = - \left(\frac{\rho_0 L^2}{\phi_0 \epsilon_t} 
         + \chi_2 \sum_i z_i c'_i\right),
     \quad \text{where } \chi_2 := \frac{e c_0 L^2}{\phi_0 \epsilon_t}.
      \label{eq.nep} 
\een
Here, the dimensionless parameter $\epsilon'$ is defined as
$\epsilon':=\epsilon/\epsilon_t$
where $\epsilon_t$ is the characteristic permittivity chosen
to be the value for water: $\epsilon_t=6.950537436\times 10^{-20}\ \text{F/\r{A}}$.
The non-dimensional parameter $\chi_2$ is approximately $125.4$ with these values. 
The corresponding dimensionless boundary conditions are
\ebe
       D'_i \left[ \frac{\partial c'_i}{\partial x'} 
          + \chi_1 \left( z_i c'_i \frac{\partial\phi'}{\partial x'} \right)
       \right] =0, 
 \quad  (\phi' - \phi'_\pm) + \eta' \frac{\partial \phi'}{\partial \nb} = 0,
  \quad \text{for } x= -1, 1,
      \label{eq.dnbc} 
\een
where $\eta' :=\eta/L$.

We drop the primes when we present our numerical methods for clarity. 


\section{Numerical Methods}
\label{secNumericalMethods}
We present a method for deriving numerical schemes 
that would conserve total concentration 
of each ion species exactly if computations were performed 
without round-off errors.  We will illustrate the method by 
describing a mass-conservative scheme (i.e. preserving ion concentration
exactly) for solving the nonlinear systems of PDEs \eqref{eq.nid} 
and \eqref{eq.nep}. The extension of the method to the multi-dimensional case 
is straightforward. This scheme uses the trapezoidal rule 
and the second-order backward differentiation formula (TR-BDF2) in time 
and the second-order central differencing in space. The TR-BDF2 scheme
is implicit in time, resulting in a system of nonlinear equations after
discretization. Instead of using the Newton-Raphson method 
for solving the large nonlinear systems at each time step, 
we present a simple iterative scheme which is easy to implement and can solve the systems efficiently.


\subsection{Discretization in Time}
\label{secDiscretizationInTime}
For time-stepping, we shall use a slight modification of the scheme described in \cite{bankcoughran1985}, which combines the trapezoidal rule with the second-order backward differentiation formula.

\noindent (1) TR step:
\begin{eqnarray}
\left\{ \begin{array}{ll} c_{i}^{n+\gamma,k+1} - \gamma \frac{\Delta t_n}{2} f(c_{i}^{n+\gamma,k+1},\phi^{n+\gamma,k}) 
&= c_{i}^n + \gamma \frac{\Delta t_n}{2} f(c_{i}^n,\phi^n),
\quad i = 1, 2, \; k =0,1,2,\dots, \\
\hfill\frac{\partial}{\partial x} \left(\epsilon
    \frac{\partial\phi^{n+\gamma,k+1}}{\partial x} \right)
&= - \left(\frac{\rho_0 L^2}{\phi_0 \epsilon_t} +\chi_2 \sum_{i=1}^2 z_i  c^{n+\gamma,k+1}_i\right), 
\label{eq.tr} 
\end{array} \right.
\end{eqnarray}
\noindent (2) BDF2 step:
\begin{align}
\left\{ \begin{array}{ll} 
c_{i}^{n+1,l+1} - \frac{1-\gamma}{2-\gamma} \Delta t_n f(c_{i}^{n+1,l+1},\phi^{n+1,l}) 
&= \frac{1}{\gamma (2-\gamma)} c_{i}^{n+\gamma}
       - \frac{(1-\gamma)^2}{\gamma (2-\gamma)} c_{i}^n, \quad 
     i = 1, 2, \; l =0,1,2, \dots, \\
\hfill\frac{\partial}{\partial x} \left(\epsilon
    \frac{\partial\phi^{n+1,l+1}}{\partial x} \right)
&= - (\frac{\rho_0 L^2}{\phi_0 \epsilon_t} +\chi_2 \sum_{i= 1}^2 z_i  c^{n+1,l+1}_i),
  \label{eq.bdf} \end{array} \right.
\end{align}
where $f(c_i,\phi)$ is defined as the right-hand side of \eqref{eq.nid}
\begin{equation}
 f(c_i,\phi) = \frac{\partial}{\partial x} 
       \left\{ D_i \left[ \frac{\partial c_i}{\partial x} 
          + \chi_1 \left( z_i c_i \frac{\partial\phi}{\partial x} \right)
       \right] \right\}.
\end{equation}
We take $\gamma = 2 - \sqrt{2}$, which minimizes 
the local truncation error.\cite{Gardner04electrodiffusionmodel}

Removing the inner iterations, corresponding to the indices $k$ 
in \eqref{eq.tr} and $l$ in \eqref{eq.bdf}, Eqs.~\eqref{eq.tr}
and \eqref{eq.bdf} is the TR-BDF2 scheme requiring a nonlinear solver
for the two systems of nonlinear equations: \eqref{eq.tr}
for $(c^{n+\gamma},\phi^{n+\gamma})$ at the grid points
and \eqref{eq.bdf} for $(c^{n+1},\phi^{n+1})$.  
With the inner iterations,  Eqs.~\eqref{eq.tr}
and \eqref{eq.bdf} provide a simple iterative scheme 
for solving the systems of nonlinear equations. For instance,
at $k$-th iteration, we update the array $c^{n+\gamma,k+1}$
at the grid points by solving the first equation
of \eqref{eq.tr} which is a tri-diagonal system after the spatial
discretization, since the values of $\phi^{n+\gamma,k}$ are known
at $k$-th iteration; then, we update $\phi^{n+\gamma,k+1}$
using the second equation of \eqref{eq.tr}. We perform
the inner iterations until convergence and, as shown later,
choosing two inner iterations $k=2$ and $l=2$ would be sufficient.
As for initial guesses at the $n$-th time step, 
we choose $\phi^{n+\gamma,0}=\phi^n$ for \eqref{eq.tr}
and $\phi^{n+1,0} = \phi^{n+\gamma,k+1}$ for \eqref{eq.bdf}
with $k$ corresponding to the last inner iteration at the previous inner
iteration.
As shall be seen in Sec.~\ref{secNumericalResults}, 
without any such inner iterations ($k=l=0$), 
one could only attain first-order accuracy in time; 
on the other hand, with just one inner iteration ($k=l=1$),
one can attain second-order accuracy in time. In other words,
the simple iterative scheme is very effective in solving the 
systems of nonlinear equations.


\subsection{Discretization in Space}
\label{secDiscretizationInSpace}
Next, we provide the discrete equations for the spatial
differential operators in Eqs.~\eqref{eq.tr} and \eqref{eq.bdf}.
Let's divide the dimensionless interval $[-1,1]$ to $J$ subintervals, 
$x_j=-1+j \Delta x$, where $\Delta x = 2/J$ and $j=0,1, \cdots, J$.
We denote the numerical values of $g(x,t)$ at $(x_j,t_n)$
by $g_j^n$ and $g(x)$ at $x_j$ by $g_j$. 
We present the standard second-order central differencing schemes
for the spatial differential operators here to facilitate
the description of the mass-conservative scheme which depends
on the details of the discretization  
at the interior grid points ($-J+1\leq j\leq J-1$).

The ion diffusion term in Eq. (\ref{eq.nid}) is discretized as
\begin{eqnarray}
\frac{\partial}{\partial x} \left(D_i \frac{\partial c_i}{\partial x}\right) (x_j) 
 \approx \frac{D_{i,j+\half} c_{j+1} - (D_{i,j+\half} + D_{i,j-\half}) c_{i,j} 
   + D_{i,j-\half} c_{i,j-1} } {(\Delta x)^2}.  
\label{eq.sdc1} 
\end{eqnarray}
The term driven by the electrostatic potential gradient
in Eq. (\ref{eq.nid}) is given by
\begin{eqnarray}   
\frac{\partial}{\partial x} \left(D_i c_i \frac{\partial \phi}{\partial x}\right) (x_j) 
\approx \frac{D_{i,j+1}c_{i,j+1} ( \phi_{j+2}-\phi_{j})-D_{i,j-1}c_{i,j-1}(\phi_{j}-\phi_{i,j-2})}{4(\Delta x)^2}.
\label{eq.scp2} 
\end{eqnarray} 
The Laplacian in the Poisson Eq.~(\ref{eq.nep}) is approximated by 
\begin{eqnarray} 
\frac{\partial}{\partial x} \left(\epsilon \frac{\partial \phi}{\partial x}\right) (x_j) 
  \approx \frac{1}{(\Delta x)^2} \left[ 
   \epsilon_{j+\half} \phi_{j+1} - (\epsilon_{j+\half} + \epsilon_{j-\half}) \phi_j 
   + \epsilon_{j-\half} \phi_{j-1}  \right].
\label{eq.sep1} 
\end{eqnarray} 


\subsection{Discretization of Boundary Condition}
\label{secDiscretizationOfBoundaryCondition}

We shall implement the boundary conditions using two different schemes. 
The first scheme is obtained by applying standard finite differencing 
to the boundary conditions, and the second is obtained by
requiring the conservation of ions within the channel. As shown later,
it is critical to preserve the ion concentrations for accurate numerical
solutions.

\subsubsection*{Standard Implementation}
\label{secFiniteDifferenceBC}
Applying the forward differencing to the right-hand side of the Nernst-Planck
equation~\eqref{eq.nid} at the left boundary 
and using the no-flux boundary condition in \eqref{eq.dnbc}, we obtain
\begin{align}
\frac{\partial}{\partial x'} 
   \left\{ D_i \left[ \frac{\partial c_i}{\partial x} 
          + \chi_1 \left( z_i c_i \frac{\partial\phi}{\partial x} \right)
       \right] \right\}(-L) 
      & \approx
       \frac{D_{i,1} \left[ \frac{c_{i,2}^{}- c_{i,0}^{}}{2\Delta x} + \chi_1z_i
       c_{i,1}^{} \frac{\phi_2-\phi_0}{2\Delta x} \right] - 0}{\Delta x} 
      \nonumber \\
      & = D_{i,1} \frac{c_{i,2}^{}- c_{i,0}^{} + \chi_1z_i
       c_{i,1}^{} (\phi_2-\phi_0)}{2 (\Delta x)^2} 
\label{eq.bc6}
\end{align}
It is similar at the right boundary. 
We implement the Robin boundary condition in \eqref{eq.dnbc} 
with the second-order central differencing using ghost grid points as
\ebe
   (\phi_0 -\phi_-) - \eta \frac{\phi_1-\phi_{-1}}{2\Delta x} = 0, 
   \quad \text{implying } \; \phi_{-1}=\phi_1 - \frac{2\Delta x}{\eta} (\phi_0-\phi_-),
\label{eq.bcep}
\een
and similarly $\displaystyle{\phi_{J+1}
=\phi_{J-1} - \frac{2\Delta x}{\eta} (\phi_J-\phi_+)}$.


\subsubsection*{Conservative Scheme: TR Step}
\label{secConservativeImplementationTR}
The no-flux boundary condition in \eqref{eq.dnbc} 
implies that the total concentration 
of each ion species is constant throughout time. 
Thus, we discretize the equations by requiring 
the numerical value of the total concentration be conserved exactly in time. 

First, we approximate the total concentration $c_{i,tot}(t_n)$ defined in 
Eq.~\eqref{eq.ctot} using the trapezoidal rule as follows
\ebe
c_{i,tot}^n= \sum_{j=1}^{J-1} c_{i,j}^n\Delta x + \frac{\Delta x}{2}\left(c_{i,0}^n+c_{i,J}^n\right)
\label{eq.ctotd}
\een
Let us examine the change of the total concentration in the TR step \eqref{eq.tr}.
\begin{eqnarray}
\frac{c_{i,tot}^{n+\gamma} - c_{i,tot}^n}{\gamma\Delta t} &=& 
    \sum_{j=1}^{J-1}  \frac{c_{i,j}^{n+\gamma} - c_{i,j}^n}{\gamma\Delta t} \Delta x +\frac{\Delta x}{2} \left( \frac{c_{i,0}^{n+\gamma} - c_{i,0}^n}{\gamma\Delta t} + \frac{c_{i,J}^{n+\gamma} - c_{i,J}^n}{\gamma\Delta t} \right)  \nonumber 
 \\
&=&\sum_{j=1}^{J-1} \left[\frac{D_{i,j+\frac{1}{2}}c_{i,j+1}^{n+\gamma}-(D_{i,j+\frac{1}{2}}+D_{i,j-\frac{1}{2}})c_{i,j}^{n+\gamma}+D_{i,j-\frac{1}{2}}c_{i,j-1}^{n+\gamma}}{2\Delta x}\right.
 \nonumber 
\\
&&+ \chi_1 z_i  \frac{D_{i,j+1}c_{i,j+1}^{n+\gamma}\left(\phi_{j+2}^n-\phi_j^n\right)-D_{i,j-1}c_{i,j-1}^{n+\gamma}\left(\phi_{j}^n-\phi_{j-2}^n\right)}{8\Delta x} 
\nonumber 
\\
 &&+ \frac{D_{i,j+\frac{1}{2}}c_{i,j+1}^{n}-(D_{i,j+\frac{1}{2}}+D_{i,j-\frac{1}{2}})c_{i,j}^{n}+D_{i,j-\frac{1}{2}}c_{i,j-1}^{n}}{2\Delta x}
 \nonumber 
\\
&&+ \left. \chi_1 z_i  \frac{D_{i,j+1}c_{i,j+1}^{n}\left(\phi_{j+2}^n-\phi_j^n\right)-D_{i,j-1}c_{i,j-1}^{n}\left(\phi_{j}^n-\phi_{j-2}^n\right)}{8\Delta x} \right]
 \nonumber
 \\
&&+ \frac{\Delta x}{2} \left( \frac{c_{i,0}^{n+\gamma} - c_{i,0}^n}{\gamma\Delta t} + \frac{c_{i,J}^{n+\gamma} - c_{i,J}^n}{\gamma\Delta t} \right) 
\end{eqnarray} 
This summation has a telescoping effect where most of the interior terms 
cancel each other and we are left with 
\begin{equation}
\begin{split}
&\frac{c_{i,tot}^{n+\gamma}-c_{i,tot}^{n}}{\gamma \Delta t}
=\frac{\Delta x}{2}\left(\frac{c_{i,0}^{n+\gamma}-c_{i,0}^n}{\gamma \Delta t}
+ \frac{c_{i,J}^{n+\gamma}-c_{i,J}^n}{\gamma \Delta t}\right) +
\\
& \frac{D_{i,\frac{1}{2}}(c_{i,0}^{n+\gamma}+c_{i,0}^{n}-c_{i,1}^{n+\gamma}-c_{i,1}^{n}) 
+ D_{i,J-\frac{1}{2}}(c_{i,J}^{n+\gamma}+c_{i,J}^{n}-c_{i,J-1}^{n+\gamma}-c_{i,J-1}^{n}) }{2 \Delta x}
\\
&-\chi_1z_i\frac{D_{i,0}(c^{n+\gamma}_{i,0}+c^{n}_{i,0})\left(\phi_{1}^n-\phi_{-1}^n\right)
+ D_{i,1}(c^{n+\gamma}_{i,1}+c^{n}_{i,1})\left(\phi_{2}^n-\phi_{0}^n\right) }{8 \Delta x}
\\
 &+\chi_1z_i\frac{D_{i,J-1}(c^{n+\gamma}_{i,J-1}+c^{n}_{i,J-1})\left(\phi_{J}^n-\phi_{J-2}^n\right)
+ D_{i,J}(c^{n+\gamma}_{i,J}+c^{n}_{i,J})\left(\phi_{J+1}^n-\phi_{J-1}^n\right)}{8 \Delta x}.
\end{split}
\label{eq.ctt}
\end{equation}
We can achieve the conservation of the total concentration 
$c_{i,tot}^{n+\gamma}=c_{i,tot}^{n}$, if we discretize 
the Nerst-Planck equation \eqref{eq.nid} at the left boundary
\begin{eqnarray}
\frac{c_{i,0}^{n+\gamma}-c_{i,0}^n}{\gamma \Delta t} 
&=& \frac{D_{i,\frac{1}{2}}(c_{i,1}^{n+\gamma} - c_{i,0}^{n+\gamma}+c_{i,1}^{n} - c_{i,0}^{n})}{(\Delta x)^2} \label{eq.trb0} \\
&&+ \chi_1z_i\frac{D_{i,0}(c^{n+\gamma}_{i,0}+c^{n}_{i,0})\left(\phi_{1}^n-\phi_{-1}^n\right)
+ D_{i,1}(c^{n+\gamma}_{i,1}+c^{n}_{i,1})\left(\phi_{2}^n-\phi_{0}^n\right)}{4(\Delta x)^2}, \hspace{.8in} \nonumber 
\end{eqnarray} 
and at the right boundary
\begin{eqnarray}
\frac{c_{i,J}^{n+\gamma}-c_{i,J}^n}{\gamma \Delta t} 
&=& -\frac{D_{i,J-\frac{1}{2}}(c_{i,J}^{n+\gamma}-c_{i,J-1}^{n+\gamma}+c_{i,J}^{n}-c_{i,J-1}^{n})}{(\Delta x)^2}\label{eq.trb1} \\
&& - \chi_1z_i\frac{D_{i,J-1}(c^{n+\gamma}_{i,J-1}+c^{n}_{i,J-1})\left(\phi_{J}^n-\phi_{J-2}^n\right)
  + D_{i,J}(c^{n+\gamma}_{i,J}+c^{n}_{i,J})\left(\phi_{J+1}^n-\phi_{J-1}^n\right)}{4(\Delta x)^2}. \nonumber 
\end{eqnarray}

It is important to point out that Eq.~\eqref{eq.trb0} 
can be seen as discretizing Eq.~\eqref{eq.nid}
using a first-order finite difference with grid size $\Delta x/2$ 
and using the no-flux boundary condition \eqref{eq.dnbc}. 
Eq.~\eqref{eq.trb0} can be rewritten as 
\begin{eqnarray}
 \frac{c_{i,0}^{n+\gamma}-c_{i,0}^n}{\gamma \Delta t} 
& = &  \frac{\left[ D_{i,\half} (c^{n+\gamma}_{i,1}- c^{n+\gamma}_{i,0})/\Delta x + 
     \frac{\chi_1z_i}{2}\left(D_{i,0}c^{n+\gamma}_{i,0}\frac{\phi_{1}^n-\phi_{-1}^n}{2\Delta x}
+ D_{i,1}c^{n+\gamma}_{i,1}\frac{\phi_{2}^n-\phi_{0}^n}{2 \Delta x}\right) \right]
    - 0   }{\Delta x}  \nonumber  \\
&& + \frac{\left[D_{i,\half} (c^{n}_{i,1}- c^{n}_{i,0})/\Delta x + 
     \frac{\chi_1z_i}{2}\left(D_{i,0}c^{n}_{i,0}\frac{\phi_{1}^n-\phi_{-1}^n}{2\Delta x}
+ D_{i,1}c^{n}_{i,1}\frac{\phi_{2}^n-\phi_{0}^n}{2 \Delta x}\right) \right] - 0
       }{\Delta x} \label{eq.bcc3} \\ 
& \approx &  
   \frac{\frac{1}{2}\left[\left(D_i \frac{\partial c_i^{n+\gamma}}{\partial x} + \chi_1 z_i D_i
  c_i^{n+\gamma} \frac{\partial \phi^n}{\partial x} \right)(x_\half)+\left(D_i \frac{\partial c_i^{n}}{\partial x} + \chi_1 z_i D_i c_i^{n} \frac{\partial \phi^n}{\partial x} \right)(x_\half) \right] - 0}{\Delta x/2}.
\nonumber 
\end{eqnarray}


\subsubsection*{Conservative Scheme: BDF2 step}
\label{secConservativeImplementationBDF2}
We can rewrite Eq.~\eqref{eq.bdf} in such a way that 
the numerical value of the derivative of the total concentration 
becomes a linear combination of the result from the TR step and the right hand side of equation \eqref{eq.nid} evaluated at the $n+1$th time step.

\vspace{-.2in}
\begin{eqnarray}
\frac{c_j^{n+1} - c_j^{n+\gamma}}{(1-\gamma) \Delta t} &=& \frac{1-\gamma}{2-\gamma} \frac{c_j^{n+\gamma} - c_j^{n}}{\gamma \Delta t} +\frac{1}{2-\gamma} f(c_j^{n+1}) 
\label{eq.linear}
\end{eqnarray}

As with the TR step, almost all of the interior terms cancel 
in a telescoping sum, and we can require the exact conservation
of the total concentration 
$\displaystyle c^{n+1}_{i,tot} = c^{n+\gamma}_{i,tot}$
in order to obtain the discretization of the Nernst-Planck equation
\eqref{eq.nid} at the boundaries 
for the BDF2 step:
\begin{eqnarray}
 \frac{c_{i,0}^{n+1} - c_{i,0}^{n+\gamma}}{(1-\gamma)\Delta t} 
&=& \frac{1-\gamma}{2-\gamma } \left( \frac{c_{i,0}^{n+\gamma} - c_{i,0}^{n}}{\gamma \Delta t}\right)
+ \frac{2}{2-\gamma} \frac{D_{i,\frac{1}{2}}(c_{i,1}^{n+1} - c_{i,0}^{n+1})}{(\Delta x)^2}
\nonumber 
\\
&&+ \frac{\chi_1z_i}{2-\gamma} \frac{D_{i,0}c^{n+1}_{i,0}\left(\phi_{1}^n-\phi_{-1}^n\right)
+ D_{i,1}c^{n+1}_{i,1}\left(\phi_{2}^n-\phi_{0}^n\right)}{2(\Delta x)^2},  \hspace{1.6in}
\label{eq.bdfb0} 
\end{eqnarray}

\begin{eqnarray}
  \frac{c_{i,J}^{n+1} - c_{i,J}^{n+\gamma}}{(1-\gamma)\Delta t}
&=&  \frac{1-\gamma}{2-\gamma} \left( \frac{c_{i,J}^{n+\gamma} - c_{i,J}^n }{\gamma \Delta t } \right)
- \frac{2}{2-\gamma} \frac{D_{i,J-\frac{1}{2}}(c_{i,J}^{n+1}-c_{i,J-1}^{n+1})}{(\Delta x)^2}
\nonumber 
\\
&&- \frac{\chi_1z_i}{2-\gamma} \frac{D_{i,J-1}c^{n+1}_{i,J-1}\left(\phi_{J}^n-\phi_{J-2}^n\right)
  + D_{i,J}c^{n+1}_{i,J}\left(\phi_{J+1}^n-\phi_{J-1}^n\right)}{2(\Delta x)^2}. \hspace{.9in}
\label{bdfb1}
\end{eqnarray}

Equation \eqref{eq.bdfb0} can be seen as discretizing only the term $f(c_{i,j}^{n+1})$ in Eq.~\eqref{eq.bdf}
using forward difference with grid size $\Delta x/2$ and using the no-flux
boundary condition in \eqref{eq.dnbc}. Eq.~\eqref{bdfb1} can be viewed 
similarly at the right boundary. 


\section{Numerical Results}
\label{secNumericalResults}

\subsection{Validation and Convergence Results}
\label{secValidationResults}
\begin{figure}[h]
\centerline{\includegraphics[width=16cm]{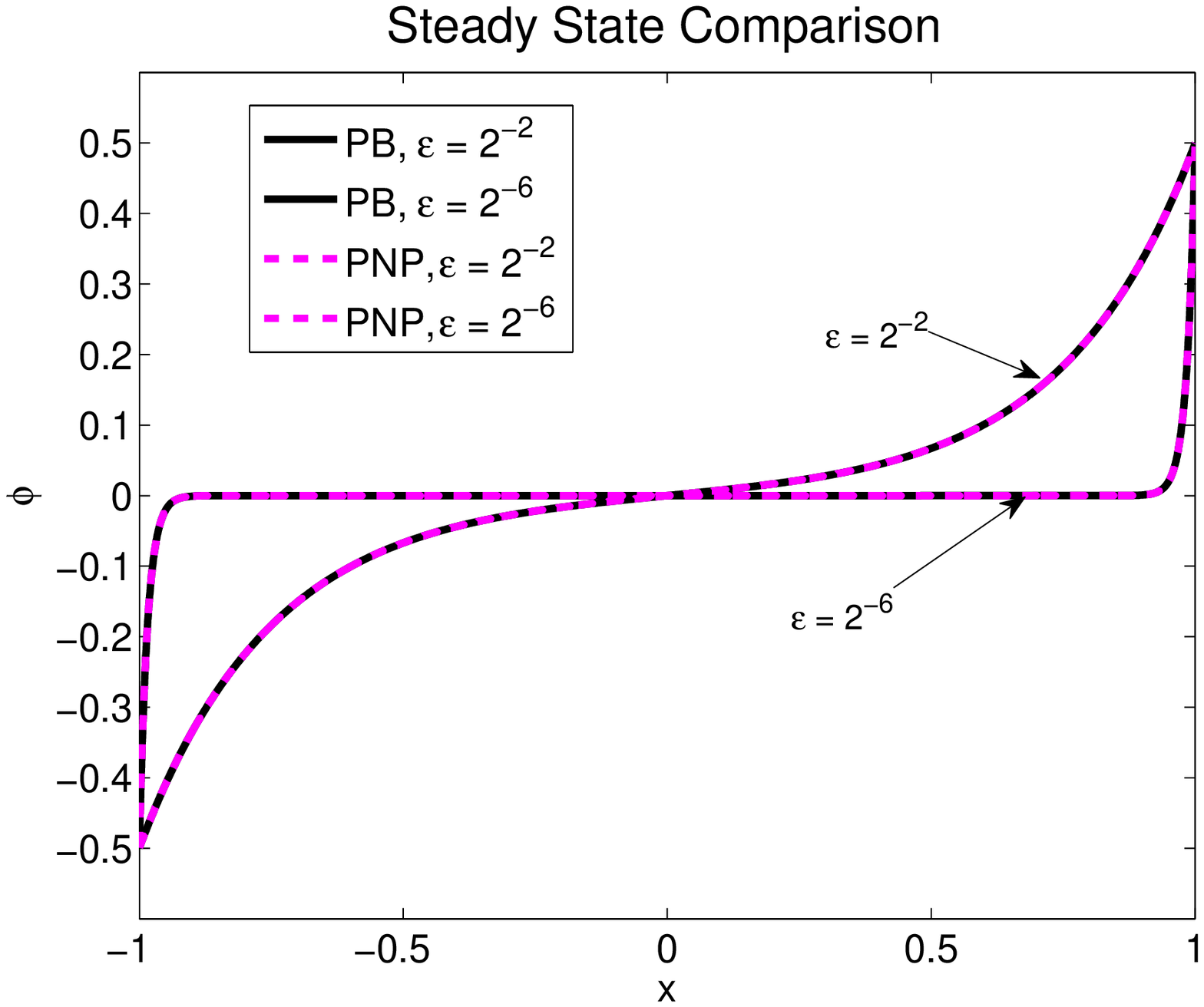}}
\caption{Comparing our steady-state solution (the dashed lines) 
using TR-BDF2 method with 
that of the Poisson-Boltzmann equation (the solid lines) obtained in \cite{liunonlinearity11}. 
The parameters are $\epsilon=2^{-2},2^{-6}$, $\eta = \epsilon$, 
$\phi_-=-1$, $\phi_+=1$. }
\label{graphs_comparison}
\end{figure}

To validate the accuracy our numerical method, we compare the steady-state 
solution from our dynamic simulations of PNP with that of the Poisson-Boltzmann
solution taken from the work \cite{liunonlinearity11}. 
Figure~\ref{graphs_comparison} shows that our steady-state solutions
match perfectly with those in \cite{liunonlinearity11} for two sets
of parameters: one with $\eta=\epsilon = 2^{-2}$ and 
the other $\eta=\epsilon = 2^{-6}$ while keeping the other parameters
constant: $\phi_-=-1, \phi_+=1, D_1=D_2=1, \chi_1=1, \chi_2=\frac{1}{2\epsilon}$, and $\rho_0=0$.  
The maximum difference in $\phi$ between the two solutions is less
than $5.6 \times 10^{-5}$. 
To get the steady-state solution, we have used the mass-conservative
TR-BDF2 method described
in previous sections with 2048 grid points in 
the interval $[-1,1]$ as in  \cite{liunonlinearity11} and 
the time-step size $10^{-4}$. At time $t=0$, the initial profiles for the ion 
concentrations are uniform in space. In this case, our time-dependent solution 
is close to the steady-state solution for the time $t\geq 2$. 
We have also verified that our solutions agree
with those in \cite{liunonlinearity11} for other sets of parameters as well, 
although they are not shown here. 

We have also checked the orders of convergence of our methods. 
The discretization method described in the previous section
 always has $O(\Delta x^2)$ convergence in space, 
regardless whether we have implemented the mass-conservative 
difference scheme or not. The order of convergence in space
is computed using the formula $\displaystyle{ 
\log_2 \frac{|\Phi(2\Delta x) - \Phi(4\Delta x)|}{|\Phi(\Delta x) - \Phi(2\Delta x)|}}$,
where $\Phi(\Delta x)$ denotes the numerical solution of the potential $\phi$
at the point $(x,t)=(0.904,0.02)$
obtained with the spatial resolution $\Delta x$. In this case, 
the time step size is chosen to be very small $\Delta t=10^{-6}$ so that 
the discretization error is dominated by that in space. 

To obtain the numerical orders of convergence in time, 
we compute the numerical solutions with three different time-step sizes
$\Delta t, 2\Delta t$ and $4\Delta t$ and then calculate
the numerical order of convergence $p$ by computing the ratio
$\Phi(2\Delta t) - \Phi(4 \Delta t)) / (\Phi(\Delta t) - \Phi(2 \Delta t))$
at the fixed position and time $(x,t)=(0.904,0.02)$. Here, the spatial
resolutions in these simulations are kept the same, $\Delta x=0.002$. 
The numerical convergence results in time are 
given in Table \ref{convresultstime}. 
We find that, if one did not perform inner iterations ($k=0$ in \eqref{eq.tr}
and $l=0$ in \eqref{eq.bdf}),
the convergence of TR-BDF2 would be first-order in time.
If we include at least one inner iteration ($k\geq 1$ and $l\geq 1$),
then the convergence becomes second-order as expected.

\begin{table}[h]
\centerline{\begin{tabular}{|c|c|c|c|}
\hline
$\Delta t$&$5\times 10^{-5}$&$2.5\times 10^{-5}$&$1.25\times 10^{-5}$\\
\hline
order of convergence for TR-BDF2, no inner loops&1.0016&1.0008&1.0028\\
order of convergence for TR-BDF2, two inner loops&2.2197&2.1779&2.2143\\
\hline
\end{tabular}}
\caption{The numerical order of convergence in time for the mass-conservative 
TR-BDF2 method solving the PNP equations in one dimension for two ion species. 
The non-dimensionalized physical parameters are $\epsilon=1$, $\eta = 4.63\times 10^{-5}$, $\phi_-=1$, $\phi_+=-1$. The calculations are performed with  $\Delta x = 0.002$
and the numerical solution of $\phi$ is evaluated at the point $(x,t)=(0.904,0.02)$.}
\label{convresultstime}
\end{table}

\subsection{Evolution of the Distributions of the Ions}
\label{secEvolution}

First, we examine the evolution of the ion concentrations and the electrostatic
 potential starting from a uniform ion distribution of two ion species
of opposite valence $z_1=1$ and $z_2=-1$: $c_i(x,0) = 1$, $i=1,2$, 
for $-1\leq x \leq 1$.  The prescribed electrostatic
potentials on the left and the right at far-field are $\phi_-=1$ 
and $\phi_+=-1$ respectively. The physical parameters
are specified as in Sec.~\ref{secNonDimensionalization}.
In the rest of this work, unless we specify otherwise, 
the non-dimensionalized parameters are chosen as
$D_1=D_2=1, \chi_1=3.1, \chi_2=125.4$ and $\eta=4.63\times 10^{-5}$,
as they were defined in Sec~\ref{secNonDimensionalization}.
Due to the symmetries of the initial and boundary conditions, the parameters 
and the domain, the profiles for the concentrations of the
two ion species at any time are symmetric with respect to the center 
of the channel, $x=0$. 

\begin{figure}[h]
\centerline{\includegraphics[width=12cm]{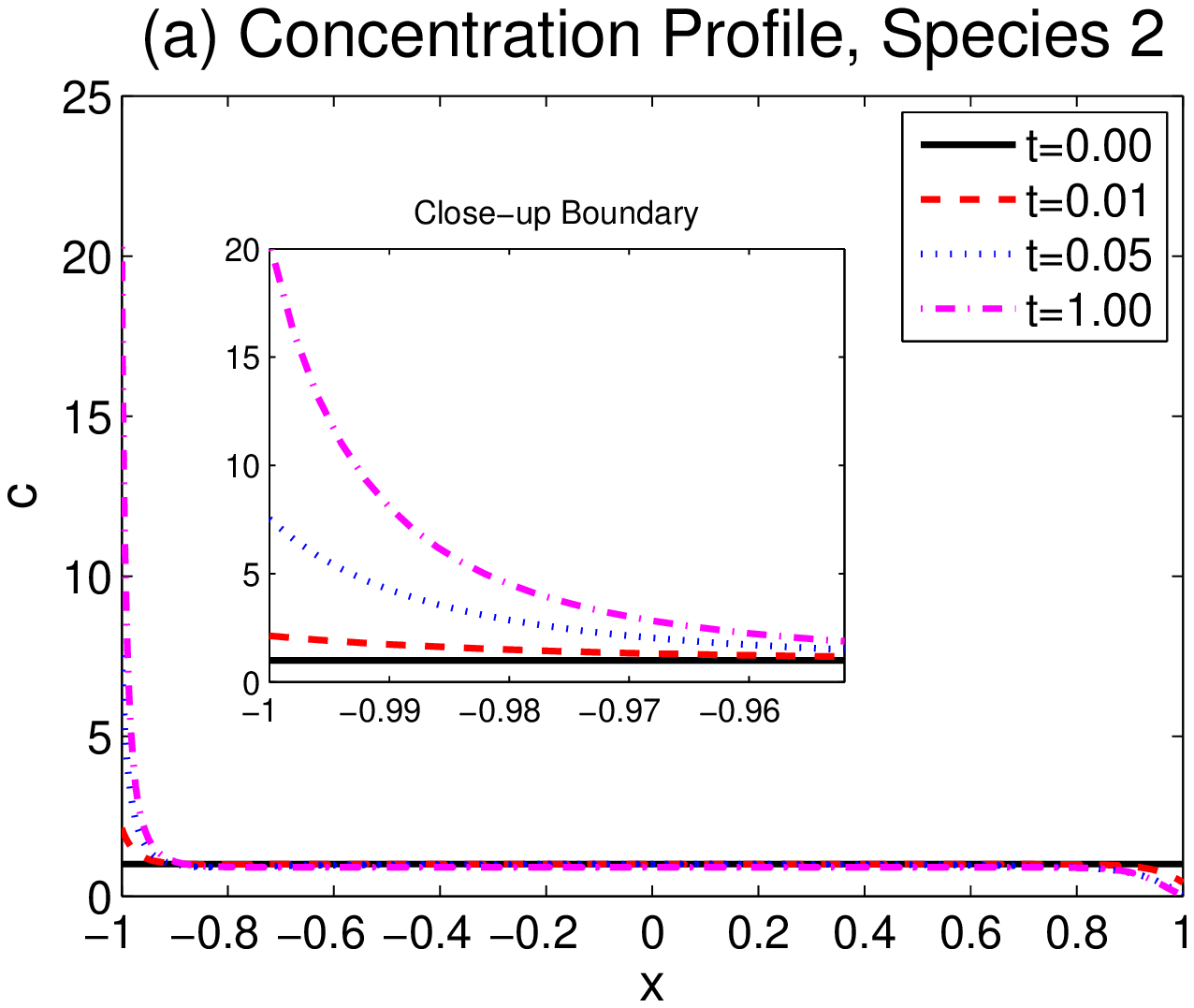}}
\centerline{\includegraphics[width=12cm]{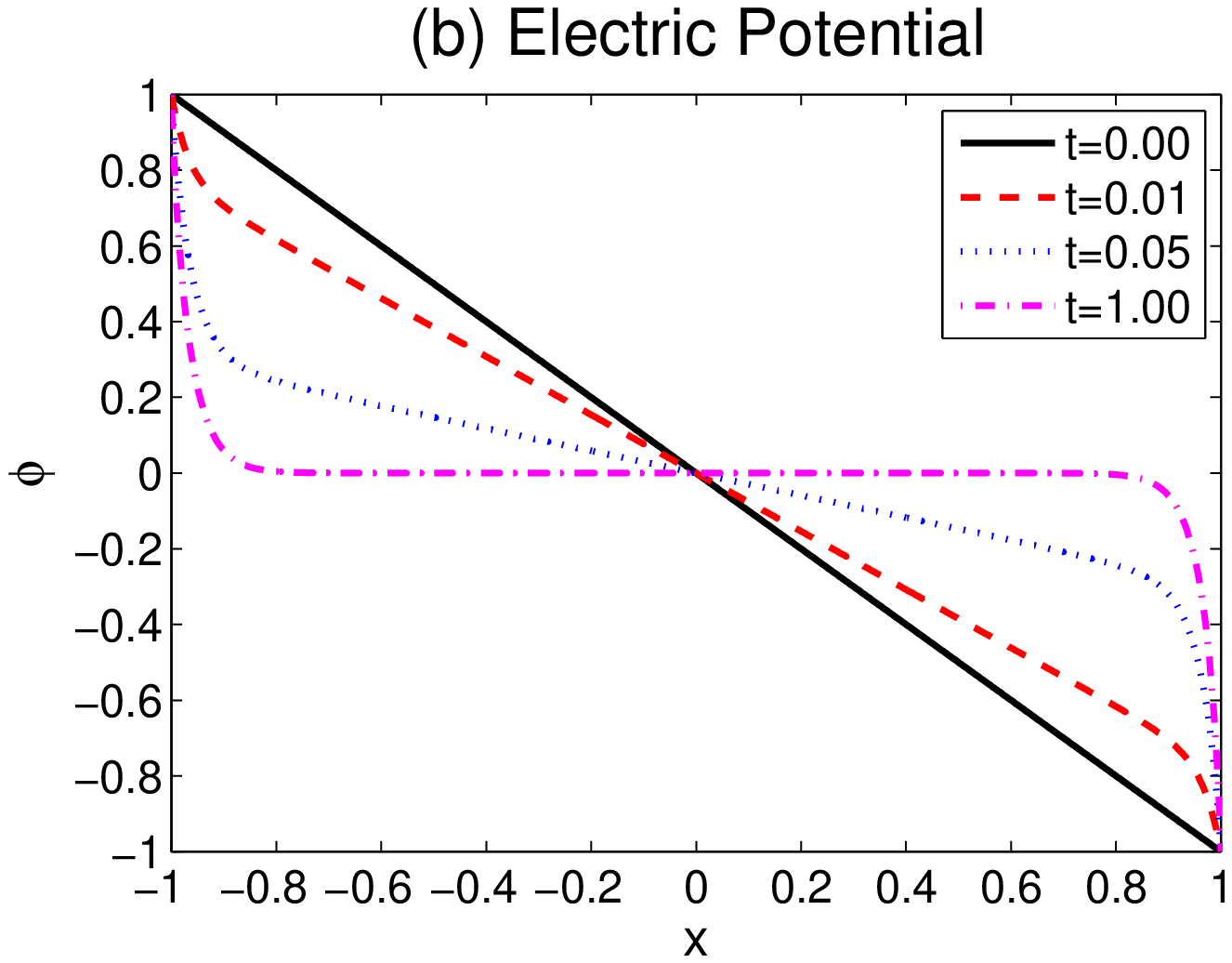}}
\caption{Simulation results using the mass-conservative TR-BDF2 method
for $\epsilon=1$, $\eta = 4.63\times 10^{-5}$, $\phi_-=1$, $\phi_+=-1$.  
The calculations were performed with $\Delta t = 10^{-4}$
and $\Delta x = 0.002$. (a) The concentration profiles for the ion species 
with the valence $z_2=-1$, $c_2(x,t)$, 
are plotted at the times $t=0$ (the solid line),
$0.01$ (dashed), $0.05$ (dotted) and $1$ (dash-dotted). 
(b) The corresponding time sequence of the electrostatic potential $\phi$
is plotted.}
\label{timelapseplots}
\end{figure}
Figure~\ref{timelapseplots} shows the profiles of the ion concentration
with the valence $z_2=-1$ and the electrostatic potential 
at the times $t=0$, $0.01$, $0.05$, and $1$. 
The Robin boundary condition \eqref{eq.dnbc} for the electrostatic
potential drives the ions with negative charges
toward the left boundary and the no-flux boundary condition \eqref{eq.dnbc} 
for the ions causes those charges to accumulate at the boundary. 
In this case, the ion concentrations keep their uniform profile 
in the bulk of the domain away from the two ends, while the electrostatic
potential changes from an initially linear profile to one that is essentially
constant (zero) except for the sharp gradient at each end. 
We find that the existence of the thin boundary layers 
requires high spatial resolution or small $\Delta x$ in the simulation. 
The numerical results would be far away from the correct solution 
if we chose $\Delta x > 0.05$.
These results show the overall behavior of the system as time elapses.

\subsection{Comparison between Mass-conservative and Standard Schemes}
\label{secComparison}

\begin{figure}[h]
\centerline{\includegraphics[width=12cm]{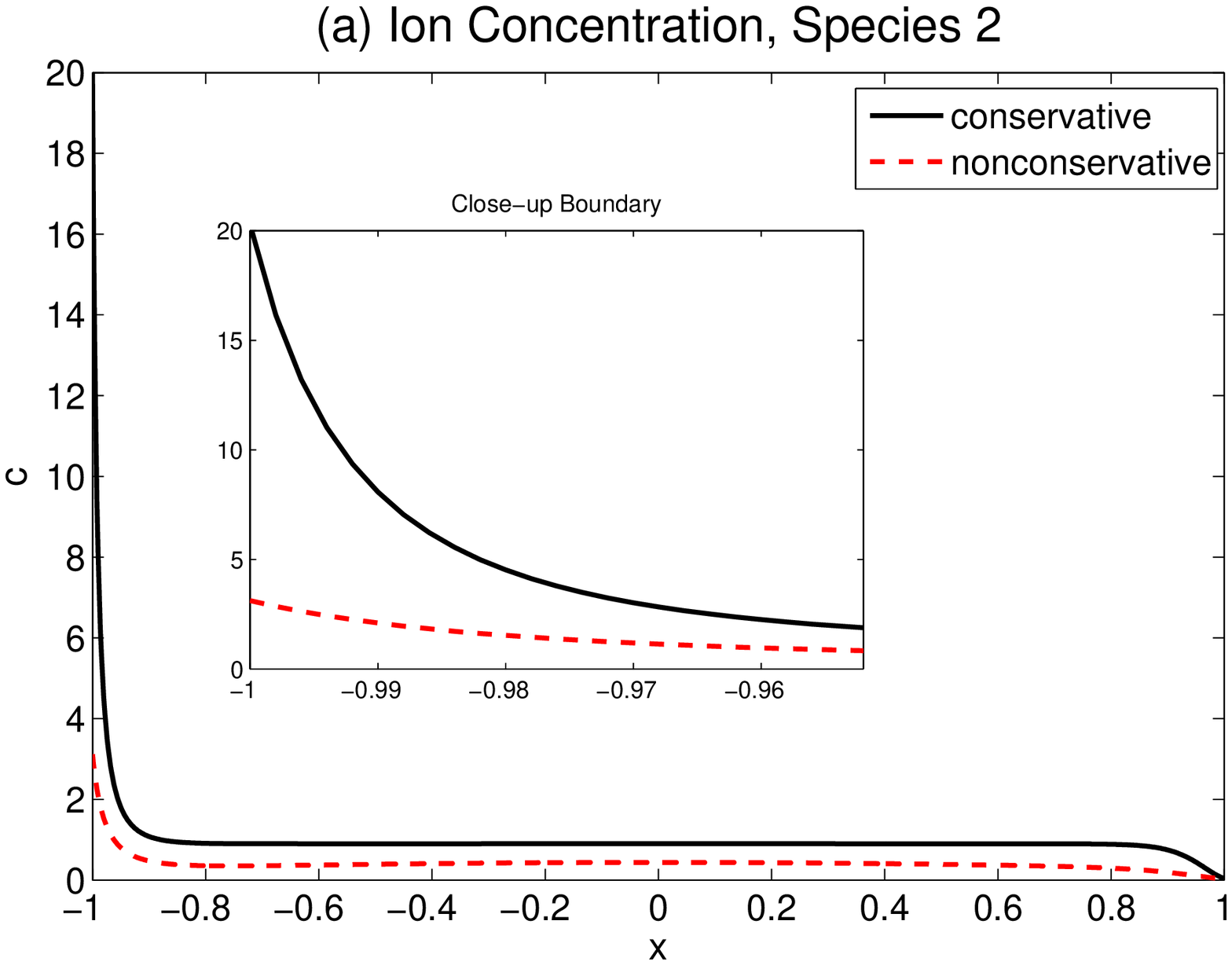}}
\centerline{\includegraphics[width=12cm]{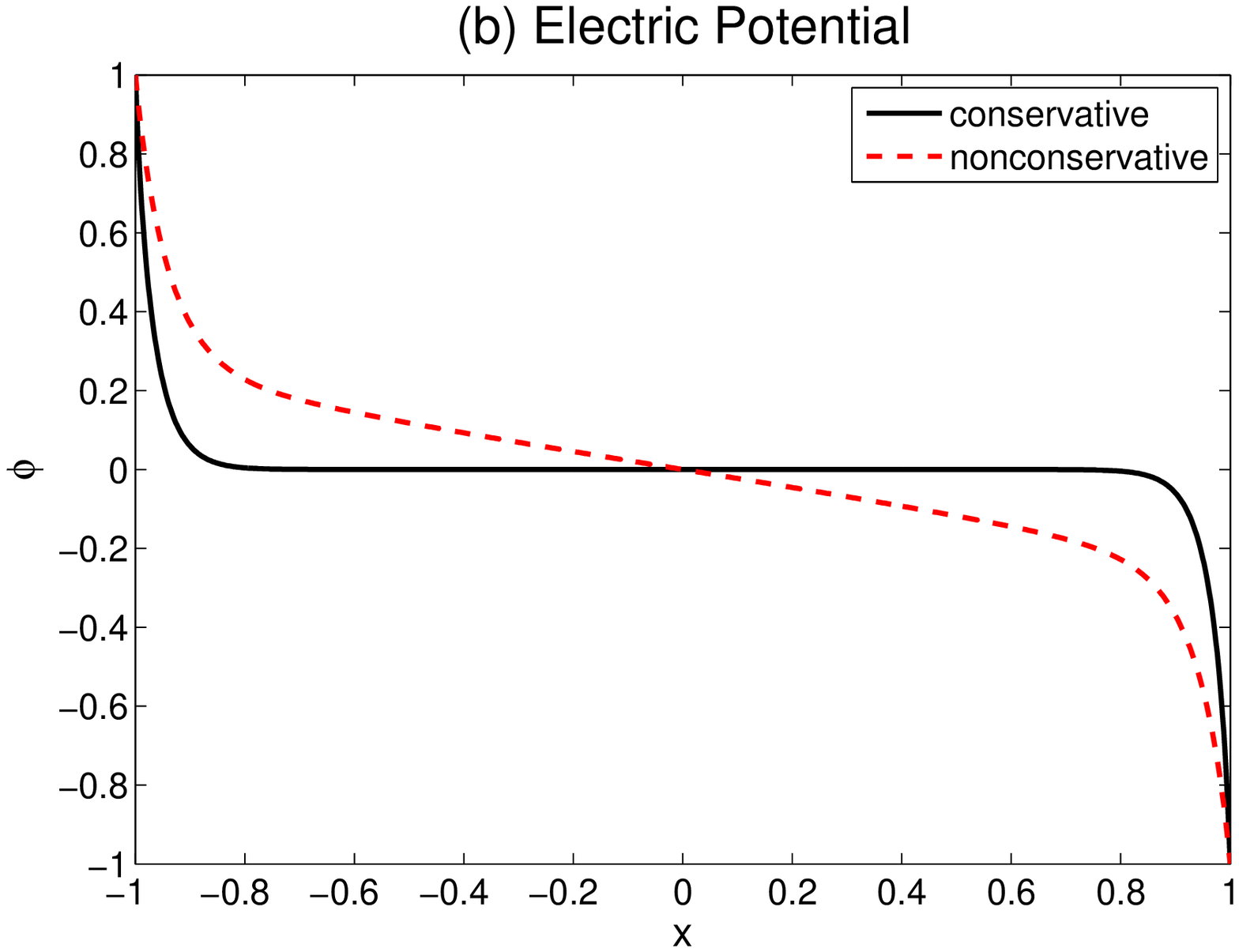}}
\caption{Comparison between the simulation results from the mass-conservative
and the non-conservative schemes for $\epsilon=1$, $\eta = 4.63\times 10^{-5}$, 
$\phi_-=1$, $\phi_+=-1$, $T=1$. The calculations were performed with 
$\Delta t = 10^{-4}$ and $\Delta x = 0.002$. 
(a) The ion concentration profiles of $c_2$ 
from the mass-conservative method (the solid line) and
the non-conservative method (the dashed line). (b) The corresponding
electrostatic potentials.}  
\label{graphs_potconc}
\end{figure}
Next, let us compare the numerical results from a standard discretization
(called as the non-conservative schemes)
of the boundary conditions, \eqref{eq.bc6}, with those
obtained from the mass-conservative schemes \eqref{eq.trb0} and 
\eqref{eq.bdfb0}. Figure~\ref{graphs_potconc} shows the ion concentration
profiles and the electrostatic potential at time $t=1$ obtained from 
both the mass-conservative schemes(the solid lines) and the non-conservative 
schemes (the dashed lines). The parameters in the computations are the same
as described in the previous Sec.~\ref{secEvolution}. To make fair
comparison, all other aspects are kept same, including 
the time-step scheme (TR-BDF2), 
the discretization scheme for interior points of the domain, 
the initial condition, the physical parameters, 
the time-step size $\Delta t$ and the space resolution $\Delta x$.   
As shown in Fig.~\ref{graphs_potconc}(a), the ion concentration 
from the non-conservative scheme is substantially lower than that 
from the mass-conservative scheme and the variations
near the boundaries are much smaller in the result 
from the non-conservative scheme.  Furthermore, the electrostatic potential
obtained from the non-conservative scheme, 
shown in Fig.~\ref{graphs_potconc}(b),
has a linear profile with non-zero slope in the middle of the domain 
and much milder slopes at the boundaries, when compared with that 
from the mass-conservative schemes.   
 
\begin{figure}[h]
\centerline{\includegraphics[width=12cm]{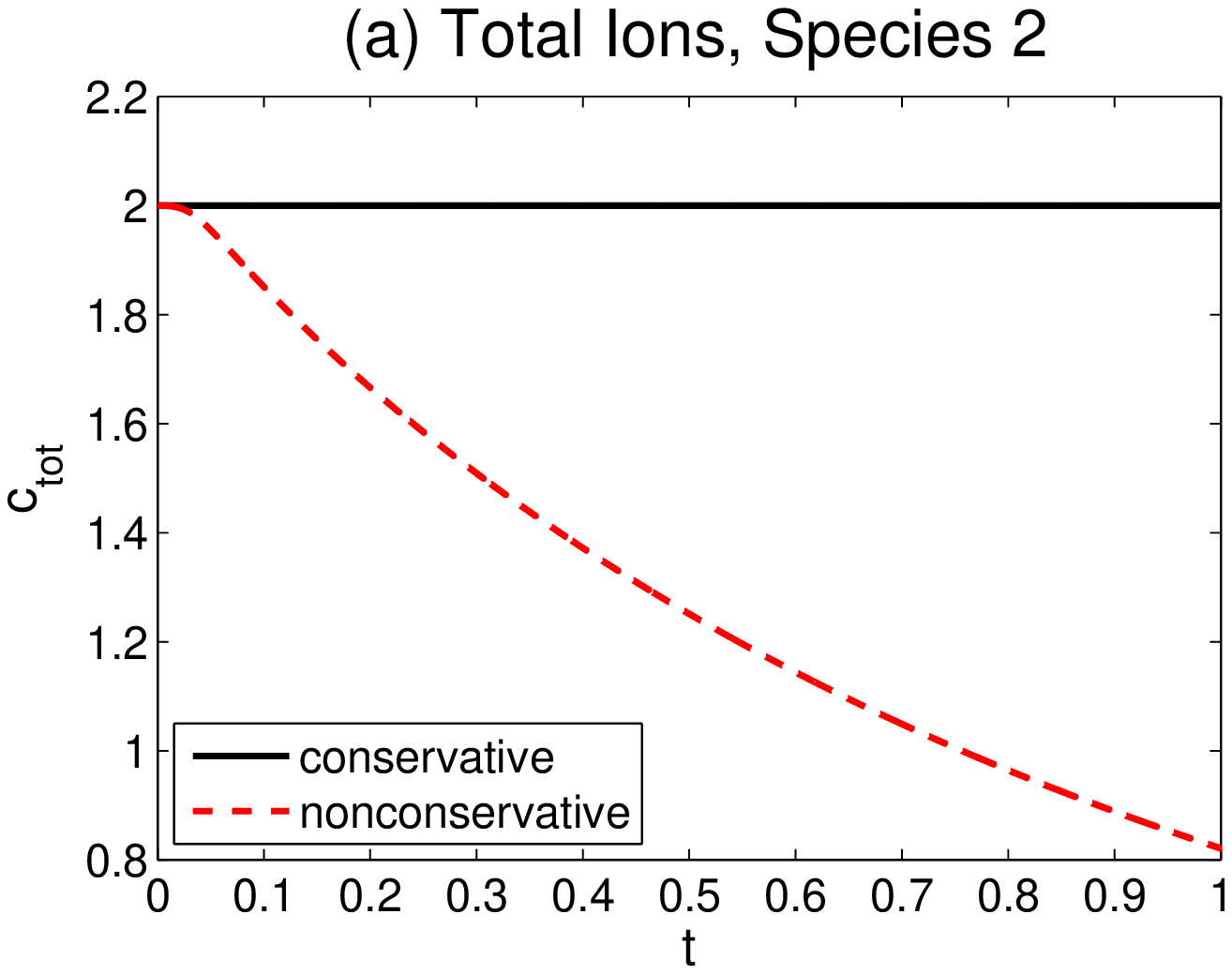}}
\centerline{\includegraphics[width=12cm]{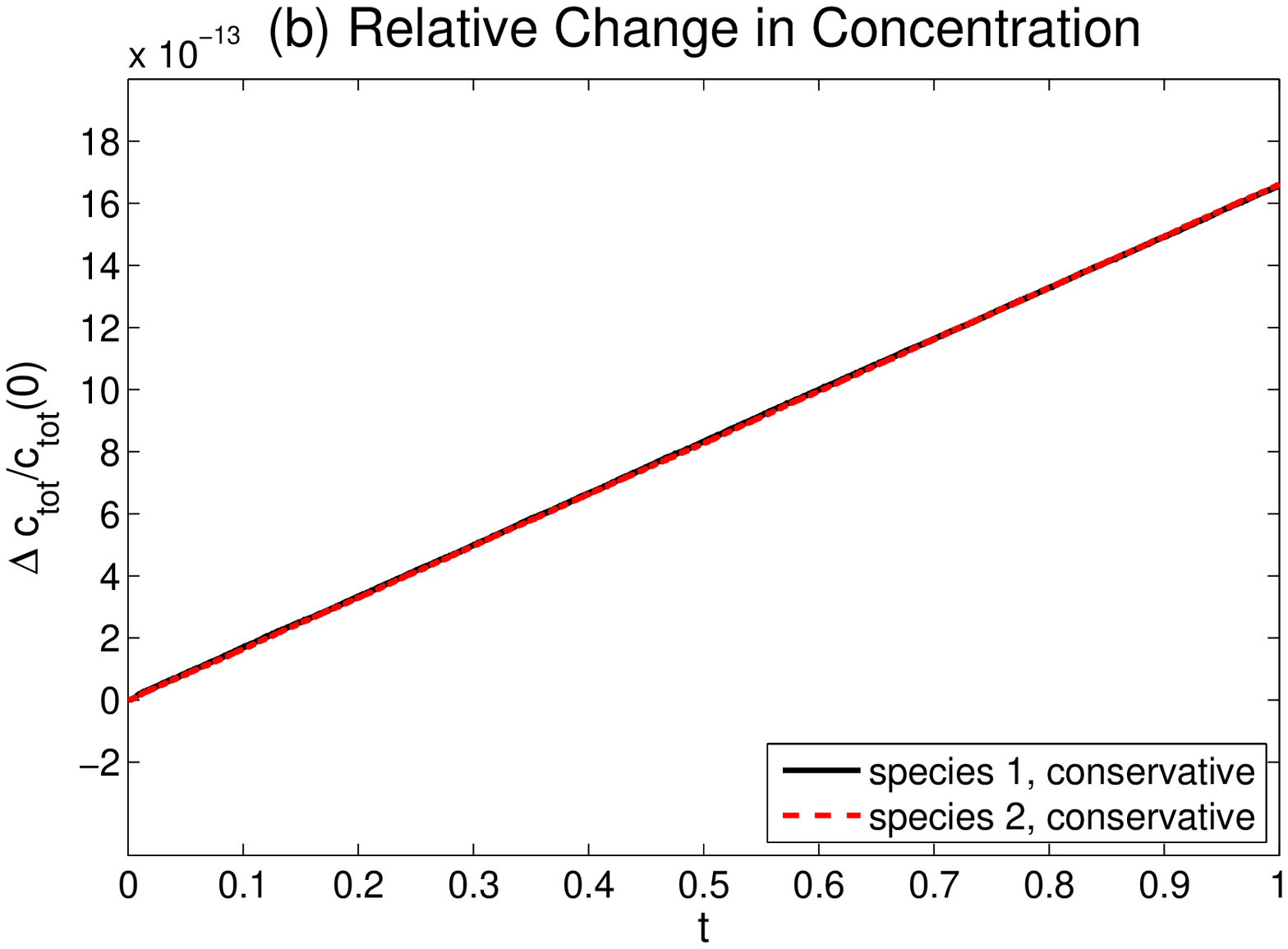}}
\caption{(a) The total ion concentration for species 2
as a function of time from the simulations using the mass-conservative (solid) 
and non-conservative (dashed) schemes. (b) The relative 
error in total concentration for both species. The parameters
are identical to those in Fig.~\ref{graphs_potconc}.} 
\label{graphs_ctot}
\end{figure}
Because of the no-flux boundary conditions \eqref{eq.npbc},
the total concentration of each ion species should be invariant in time.
Figure~\ref{graphs_ctot} shows that the mass-conservative scheme
preserves the conservation of the ions perfectly (up to the level
of roundoff error) over a long period of time, 
while the total number of ions at the time $t=1$ 
obtained from the non-conservative scheme is 
reduced to less than half of the original amount. 

\begin{figure}[h]
\centerline{\includegraphics[width=12cm]{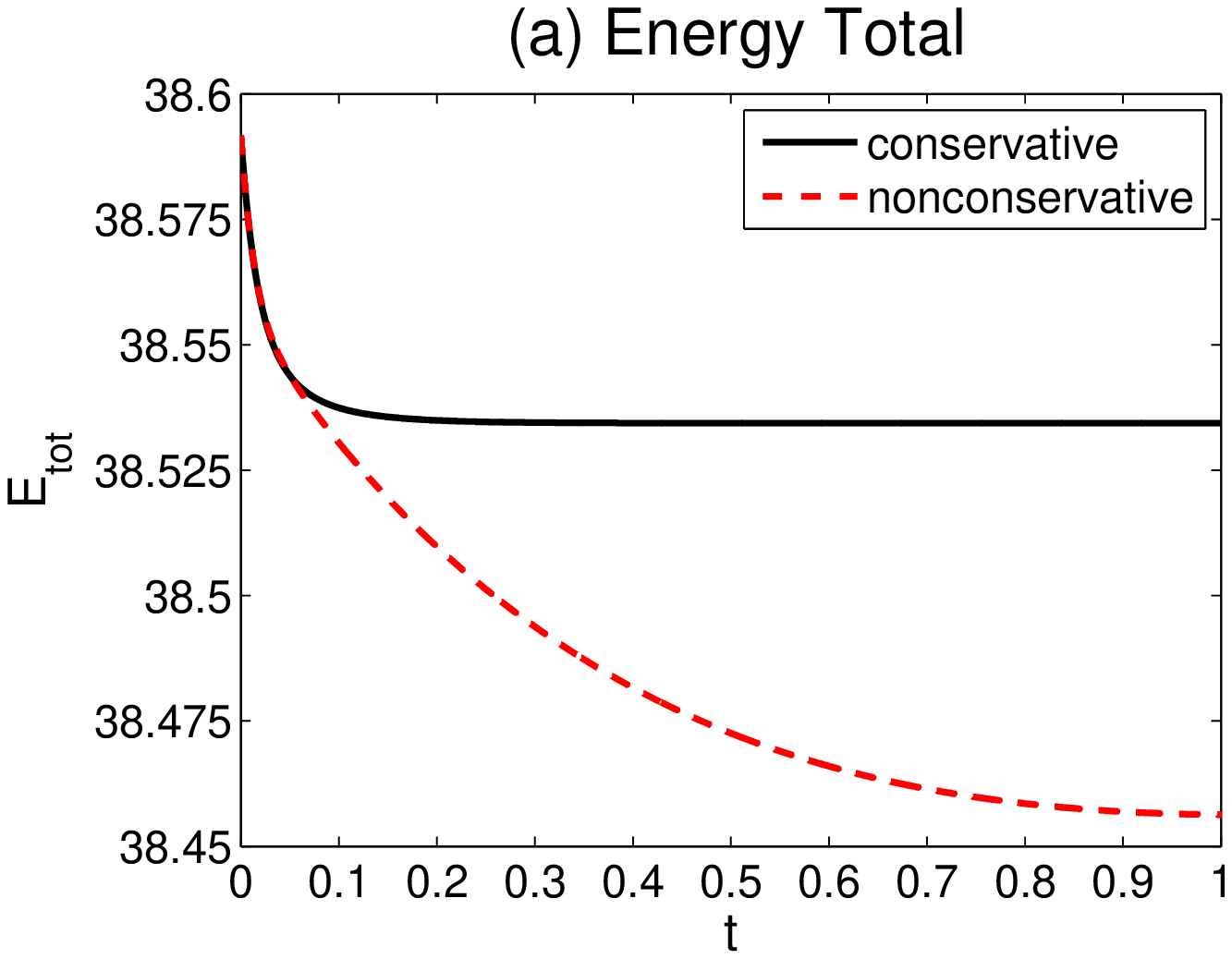}}
\centerline{\includegraphics[width=12cm]{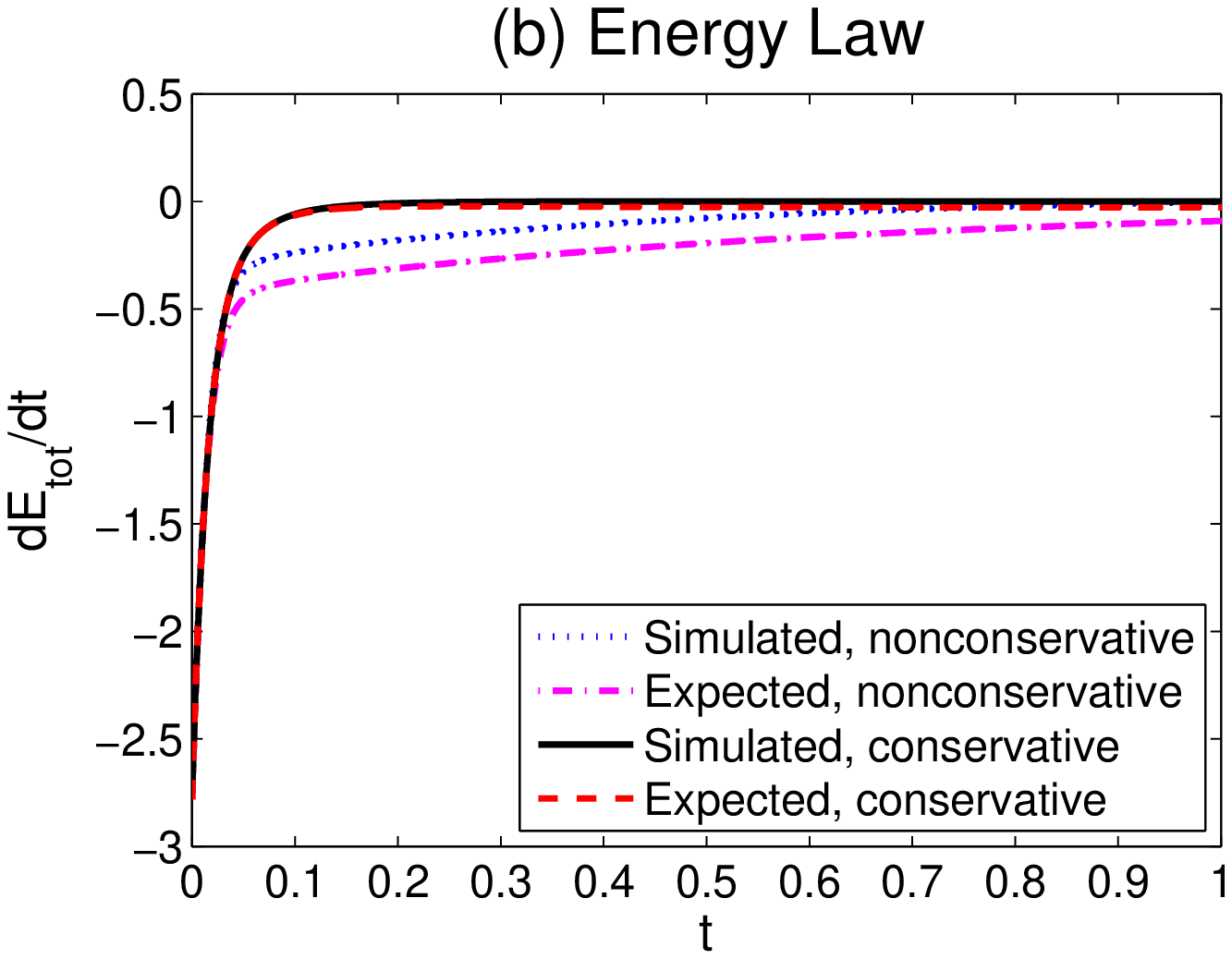}}
\caption{(a) The total energy 
as a function of time from the simulations using the mass-conservative (solid) 
and non-conservative (dashed) schemes. (b) The rate of change in energy,
$\frac{{\rm d} E}{{\rm d}t}$, obtained from the graph (a) 
and the right-hand side of Eq.~\eqref{eq.ourenergylaw}. 
The solid and the dotted lines correspond to the left-hand side
of Eq.~\eqref{eq.ourenergylaw} for the mass-conservative 
and the non-conservative schemes respectively. 
The dashed and the dash-dotted lines correspond to the right-hand side
of Eq.~\eqref{eq.ourenergylaw} for the mass-conservative 
and the non-conservative schemes respectively. 
The parameters are identical to those in Fig.~\ref{graphs_potconc}.} 
\label{graphs_energy}
\end{figure}
Figure~\ref{graphs_energy}(a) shows that the total energy $E$  
as a function of time $t$ for both the conservative and non-conservative schemes.
The total energy obtained from the mass-conservative scheme approaches
the minimum energy state much faster than that from the non-conservative
scheme.  
More importantly, in Sec.~\ref{secEnergyDissipation}, it is shown that
the total energy of the system $E$ defined as \eqref{eq.ten} satisfies
the energy dissipation law \eqref{eq.ourenergylaw}. 
In Fig.~\ref{graphs_energy}(b), we plot the rate of change in energy,
$\frac{{\rm d} E}{{\rm d}t}$, for the mass-conservative (the solid line)
and the non-conservative schemes (the dotted line) 
obtained by using a second-order finite difference
based on the numerical result $E(t)$ shown in Figure~\ref{graphs_energy}(a). 
In the same graph, we also plot the expected dissipation rate 
given by the right-hand side of \eqref{eq.ourenergylaw}, 
computed using the second-order central differencing and trapezoidal rule
and shown by the dashed line for the conservative scheme
and the dash-dotted line for the non-conservative scheme
in Fig.~\ref{graphs_energy}(b).  
It shows that the numerical result from
the conservative scheme (the solid line) agrees with the energy dissipation
law (the dashed line) very well. 
In contrast, the corresponding
results for the non-conservative scheme show that 
the energy dissipation law is not satisfied after a short period of time. 
This is due to the fact that the total concentration from the non-conservative scheme 
displays very poor performance in conserving the total concentrations. 
The results show that the discretization 
of the boundary conditions have profound impact on satisfying 
the physical properties: the energy dissipation law 
and the conservation of the total number of ions.  

\begin{figure}[h]
\centerline{\includegraphics[width=16cm]{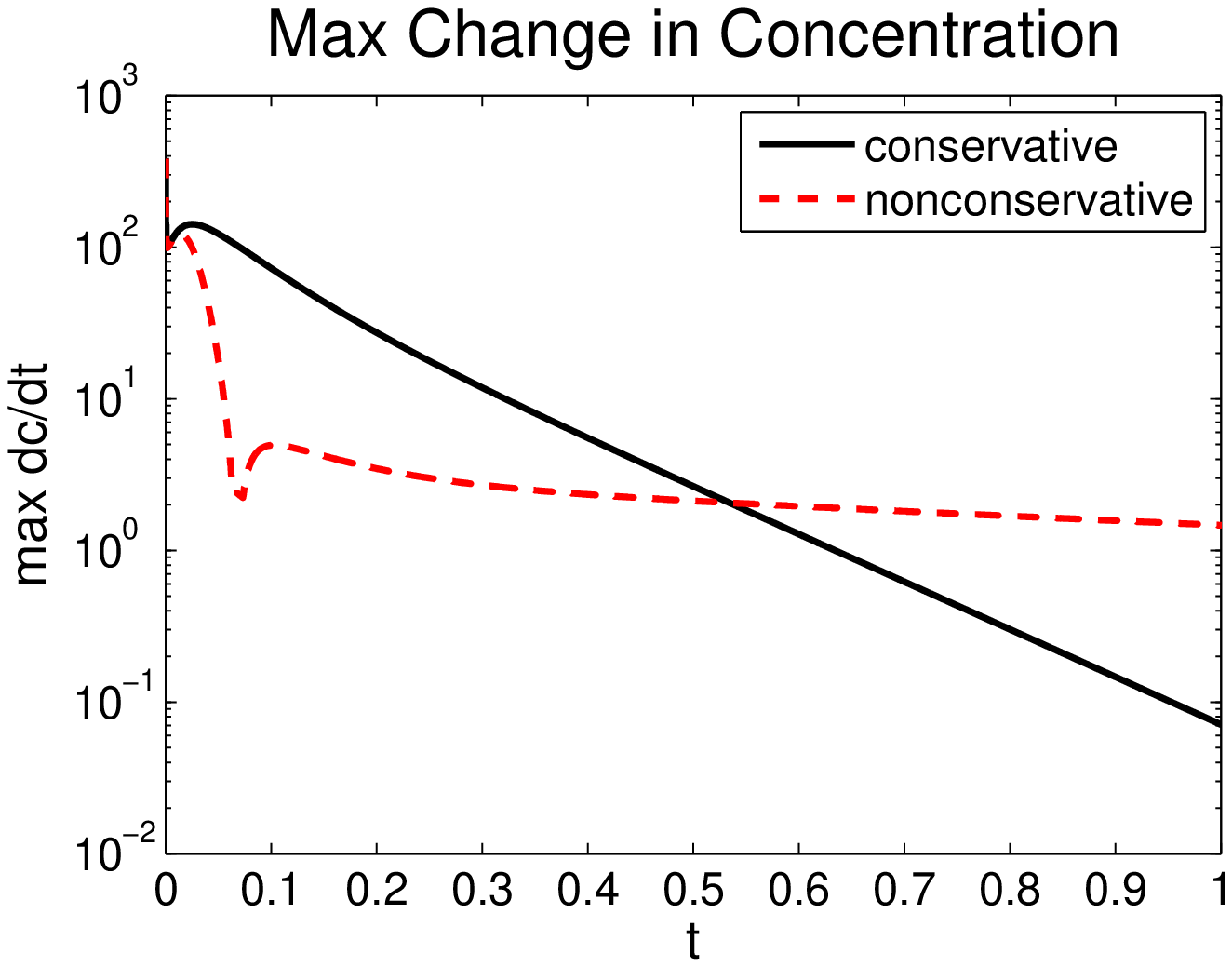}}
\caption{The maximum rate of change in ion concentrations 
as a function of time for the non-conservative(the dashed line) and 
conservative(the solid line) schemes. The parameters are identical to those in Fig.~\ref{graphs_potconc}.}
\label{graphs_maxdc}
\end{figure}
In addition to energy decay, we compute the maximum rate of change 
in the concentrations of the species over the domain, 
i.e. $\displaystyle{\max_{i, -1\leq x \leq 1} |\frac{\partial c_i}{\partial t}}|$.
It is notable from the time derivative of concentration shown 
in Fig.~\ref{graphs_maxdc} that  
the numerical results from the conservative numerical scheme
steadily approach the equilibrium in time. On the other hand, 
the non-conservative scheme is approaching a steady state much faster initially,
but, later in time, the non-conservative scheme's behavior changes 
and it does not appear to reach a steady state. This result emphasizes 
the necessity of the conservative numerical scheme for long-time simulation.

\subsection{Effect of Parameters}
\begin{figure}[h]
\centerline{\includegraphics[width=8cm]{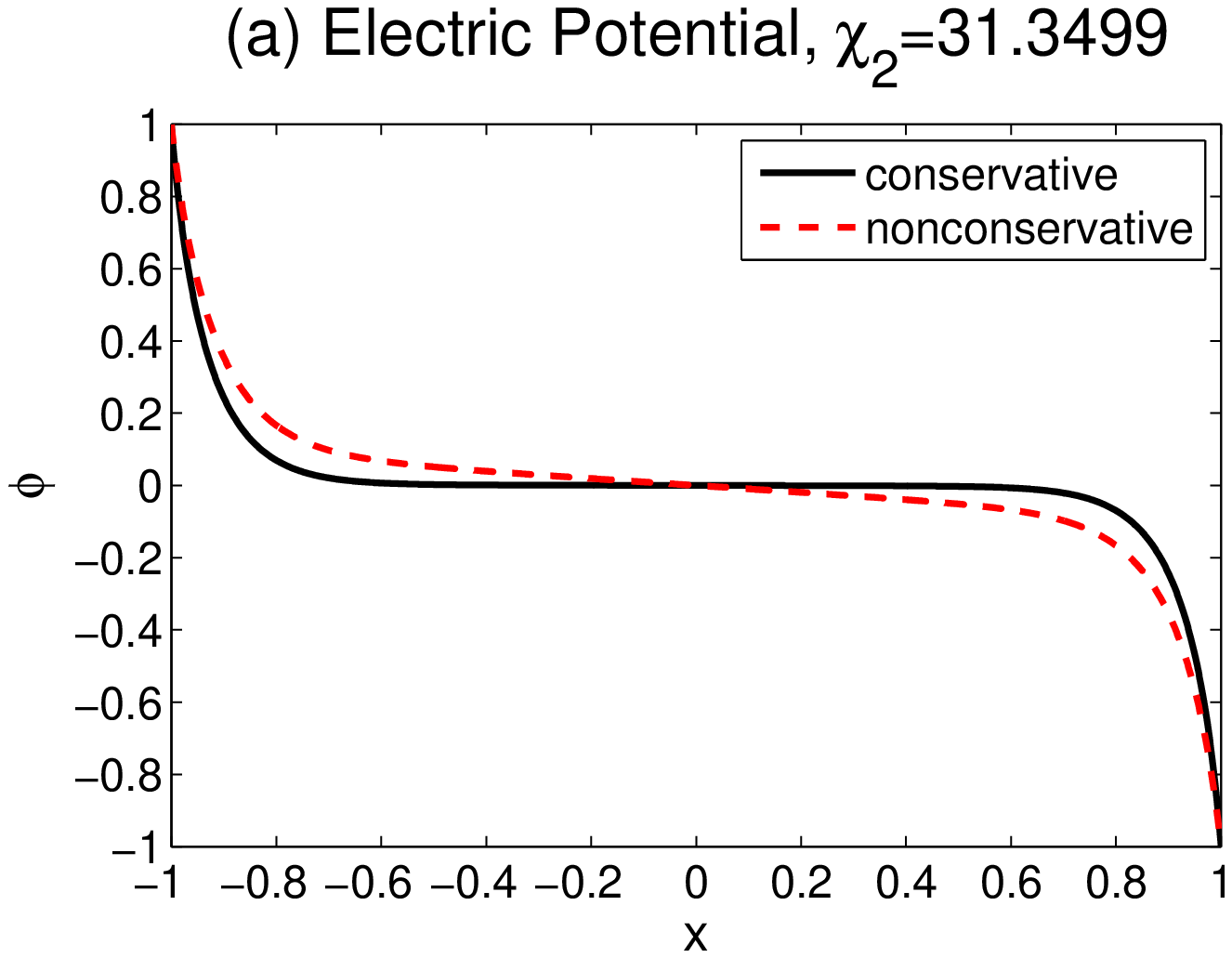}
\includegraphics[width=8cm]{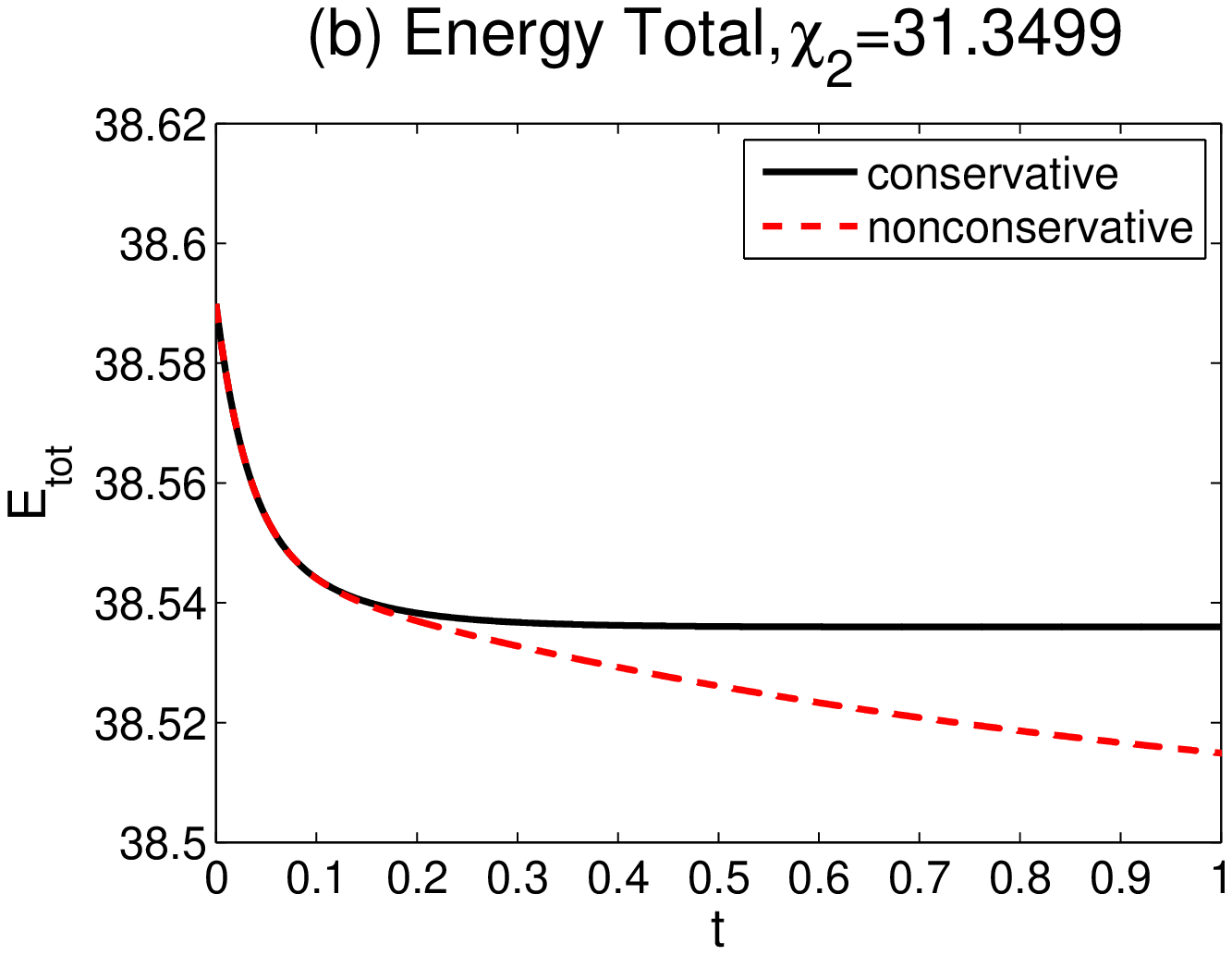}}
\centerline{\includegraphics[width=8cm]{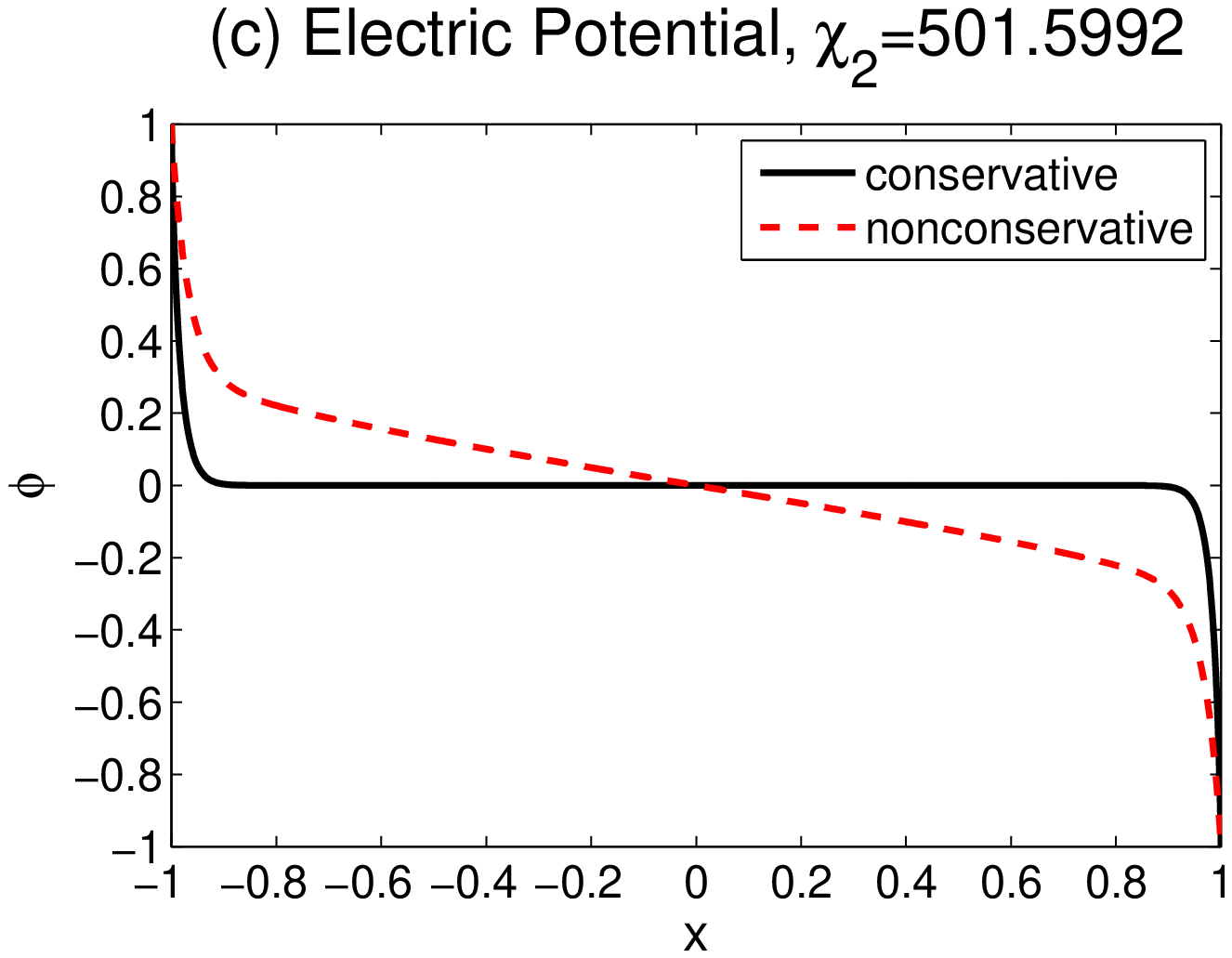}
\includegraphics[width=8cm]{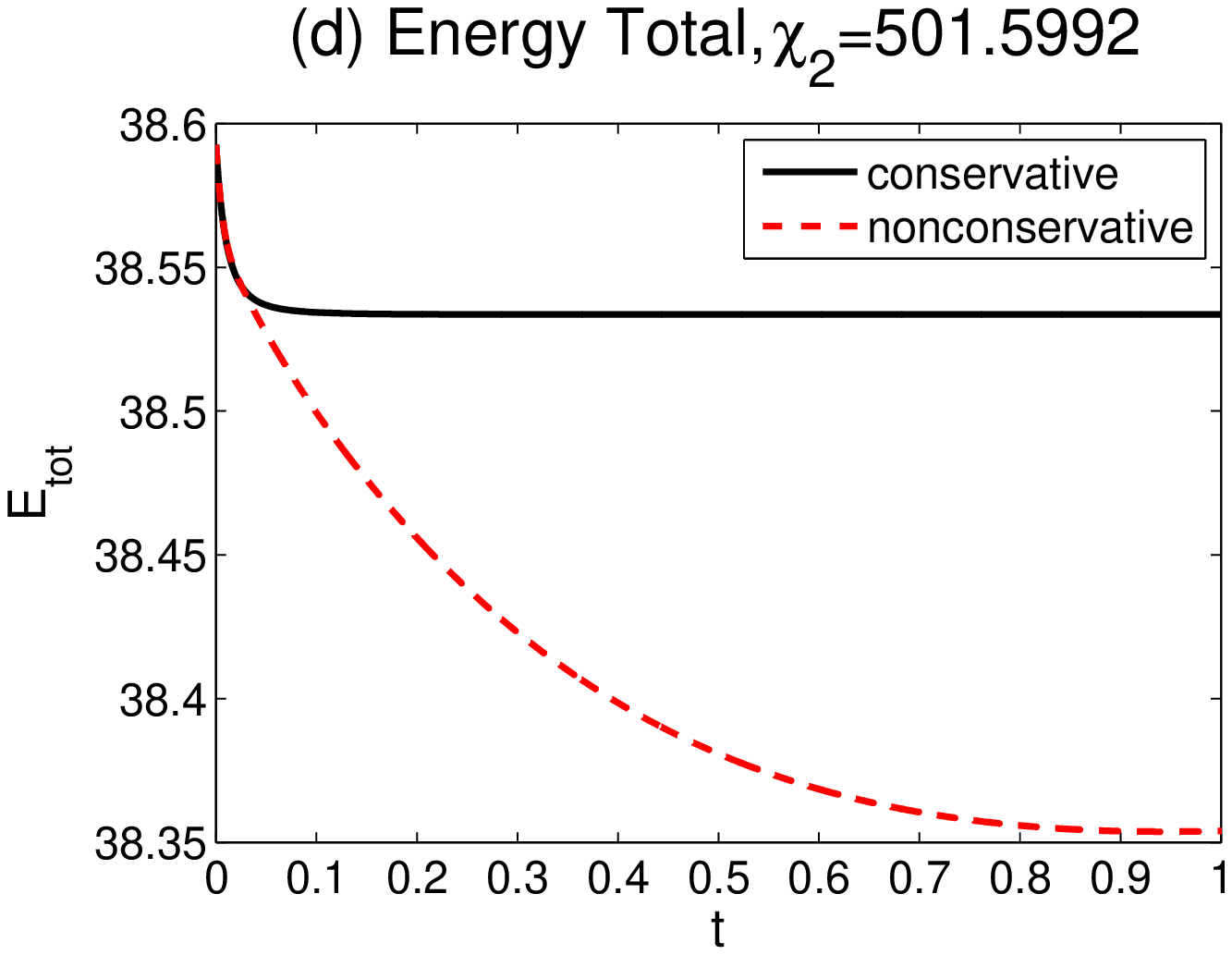}}
\caption{Comparison between the simulation results from the mass-conservative
and the non-conservative schemes for different values of the non-dimensional parameter $\chi_2$. The calculations were performed with $\Delta t = 10^{-4}$ and $\Delta x = 0.001$. The other parameters are identical to those in Fig.~\ref{graphs_potconc}.}
\label{graphs_chi2comp}
\end{figure}

The size of the difference in the results 
from conservative and non-conservative schemes 
depends on the non-dimensional parameter 
$\displaystyle \chi_2=\frac{e c_0 L^2}{\phi_0 \epsilon_t}$. 
For the physical model of the ion transportations, 
the value of $\chi_2$ can be arbitrarily large, 
depending on the values of average ion concentration $c_0$
and the applied electrostatic potential $\phi_0$ at the boundaries.
Consequently, it is important to pay attention to the size
of the dimensionless parameter $\chi_2$.
In Fig.~\ref{graphs_potconc}, we have shown that, for $\chi_2= 125.4$, 
the results of non-conservative schemes are far away from the correct 
results. Figure~\ref{graphs_chi2comp}(a) and (c) show the profiles of 
the electrostatic potential $\phi$ at a fixed time $t=1$ 
from both the conservative and the non-conservative schemes
with two more different values of $\chi_2=31.35$ and $501.6$,
while keeping all other parameters the same as those 
for Fig.~\ref{graphs_potconc}. At $t=1$, the system has reached 
the steady state, shown by the constant values  
for the conservative scheme in the energy plots
of Fig.~\ref{graphs_chi2comp}(b) and (d). Comparing the graphs of potential in 
Fig.~\ref{graphs_potconc}(a), (c) and Fig.~\ref{graphs_chi2comp},
we find that the value of $\chi_2$ primarily affects 
the width of the boundary layer, 
with larger $\chi_2$ resulting in thinner boundary layers. 
A thinner boundary layer transitions much more sharply near the boundaries, 
and thus requires more computational grid points in the region
and more truthful discretization of the boundary conditions.
This causes the differences in electrostatic potential profiles 
and the energy dissipation in time (shown by Figs.~\ref{graphs_chi2comp}(b)
and (d)) between the conservative
and non-conservative schemes to be greater as one increases $\chi_2$. 
A thinner boundary
layer also affects performance with regard to the energy dissipation law,
which is not shown here in plots. 
Larger $\chi_2$ leads to a larger discrepancy between the decay rate
of the total energy (the left-hand side of Eq.~\ref{eq.ourenergylaw})
and the energy dissipation rate (the right-hand side of the law Eq.~\ref{eq.ourenergylaw}),
and this discrepancy gets worse faster for the non-conservative scheme 
than for the conservative scheme.

\begin{figure}[h]
\begin{center} {\includegraphics[width=8.5cm]{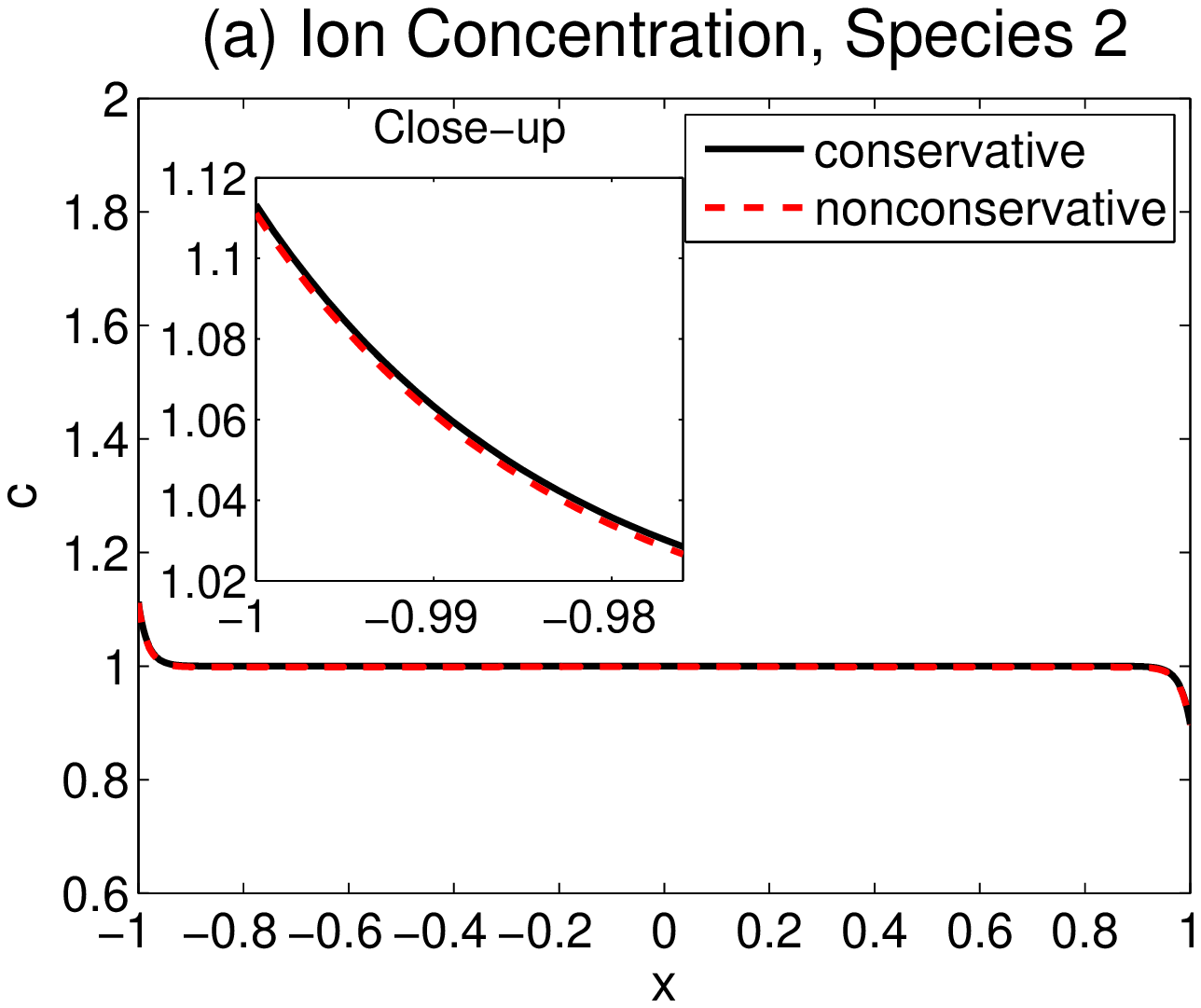}
\includegraphics[width=8.5cm]{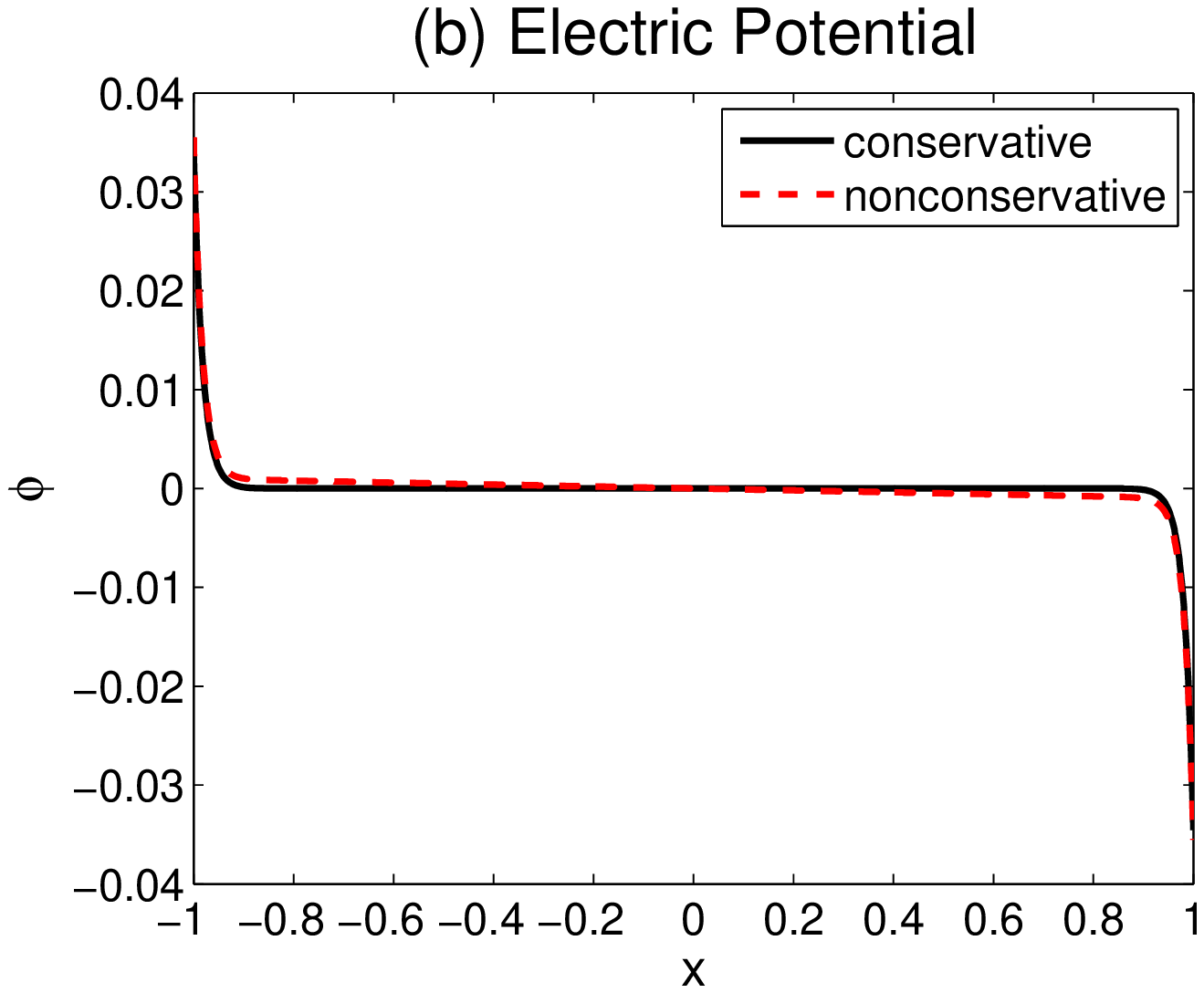}
\includegraphics[width=8.5cm]{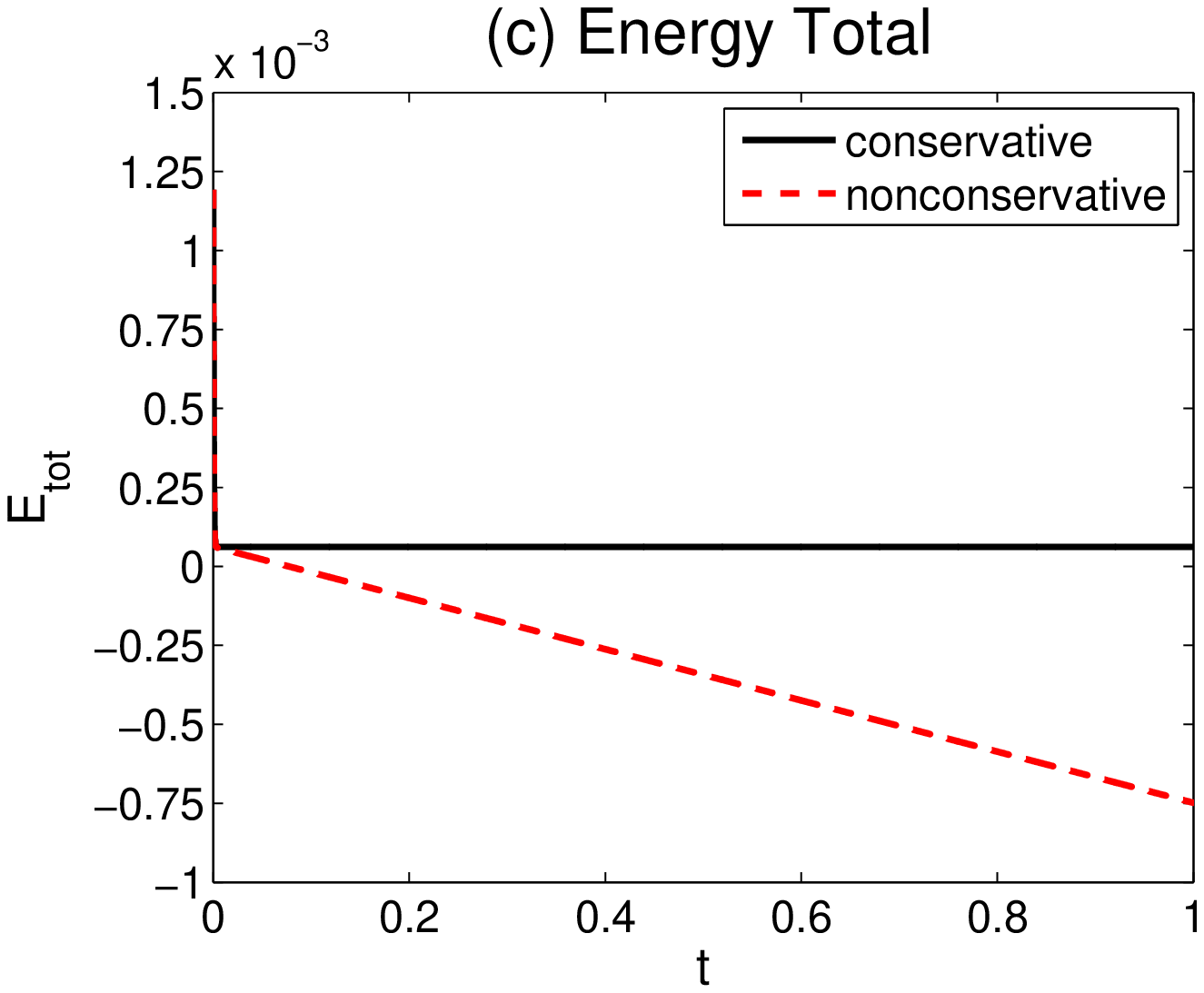}}
\end{center}
\caption{Comparison between the simulation results from the mass-conservative
and the non-conservative schemes for $\eta = 0.5$. The other parameters are identical to those in Fig. \ref{graphs_potconc}. (a) The ion concentration
at the non-dimensionalized time $T=1$. (b) The electric potential at the non-dimensionalized time $T=1$.
(c) The change of the total energy in time.}
\label{graphs_higheta}
\end{figure}

Finally, we examine the effect of the parameter $\eta$ in the Robin boundary
condition \eqref{eq.pbc} on the numerical results.
As noted in Sec.~\ref{secValidationResults}, the steady state changes 
dramatically if the relative values of the physical parameters $\eta$ 
and $\epsilon$ are changed. In order to determine the effect of $\eta$ itself 
on the results, we have tested a range of non-dimensionalized values 
for $\eta$ ranging from $10^{-6}$ to $1$, while holding $\epsilon$ at its constant non-dimensionalized value of $1$. We find that, when $\eta$ increases
from $10^{-6}$ to $0.001$, the concentration profiles at the steady state 
do not change much, having a maximum relative difference of only $10^{-4}$, but this property does not generalized to larger $\eta$. 
We also find that the discretization error, 
especially for the non-conservative scheme, is significantly 
affected by the value of $\eta$. For large values of $\eta$, say $\eta > 0.1$, 
the growth of the discretization error of the non-conservative scheme is 
rather slow, and consequently the concentration and 
electric potential profiles obtained from the non-conservative
scheme are close to those obtained 
by the mass-conservative schemes. An example of this property 
is shown in Fig.~\ref{graphs_higheta}. It appears 
that, for $\eta=0.5$, the total energy from the non-conservative scheme
decreases linearly in time after an initial sharp drop, 
becoming negative at later time. On the other hand, the conservative scheme reaches a steady state very quickly and does not deviate from it. For small values
of $\eta$ such as those shown in Fig.~\ref{graphs_potconc}, 
both the conservation property of the total concentrations 
and the energy dissipation law deteriorate at a fast pace for the non-conservative scheme, 
and the difference between the results from the conservative 
and the non-conservative schemes grows bigger as $\eta$ gets smaller. 

\section{Conclusion} \label{secConclusion}

The primary objective of this work is to investigate the effects of 
conservation property of discretization schemes on the numerical results.
We have shown that, with regard to the PNP equations, whether a numerical method
preserves the mass conservation could have a critical impact on 
the behavior of the system, especially the steady state results. 
We have provided a discretization scheme that preserves 
the mass conservation exactly (excluding the round-off errors) 
and the energy dissipation law well for long-time simulation. 

Our method is implicit in time and second-order accurate in both space and time.
We have verified that approximating the fully implicit solution is 
necessary for second-order convergence in time. Further, we find that
one can avoid using Newton-type nonlinear solvers by performing a simple iterative scheme.

In this work, we have simulated the equations with realistic 
physical parameters, particularly investigating the effect of 
the non-dimensional parameters $\chi_2$ in the Poisson equation 
and $\eta$ in the Robin boundary condition for the electrostatic
potential. We find that the mass-conserving scheme is more robust to changes in parameters,
especially changes to the value of $\eta$. 

Although this work makes good progress in constructing an accurate method 
for solving the Poisson-Nernst-Planck equations numerically, 
there are many challenges remaining. First, one of them is to account 
for the finite size of the ions as its effect is enormous considering
the narrow width of the ion channels.\cite{liu10a,Horng2012} Second, for most ion channels, 
the appropriate boundary conditions are Dirichlet-type. We will investigate
the possibility to preserve the energy dissipation law exactly instead
of the mass and study the effect of the conservation on long-term behavior
of the simulation. Third, we would like to include distributions of
permanent charges for studying selectivity of ion channels. 

\section{Acknowledgement} 
X. Li is partially supported by the NSF grant DMS-0914923 and C. Liu is partially supported by the NSF grants 
DMS-1109107, DMS-1216938 and DMS-1159937.

\bibliography{ic}

\end{document}